\documentclass[11pt,a4paper,Color]{article}
\usepackage{amssymb, paralist}
\usepackage{amsfonts, amsmath, wasysym}
\usepackage{latexsym}
\usepackage{graphicx, color}
\usepackage{indentfirst}

\usepackage{amssymb}
\usepackage{amsmath}
\usepackage{graphicx, color}
\usepackage{indentfirst}

\def\qed{\hfill {\large ${\sqcup\!\!\!\!\sqcap}$}}

\newenvironment{proof}{{\bf Proof. }}
{\qed \\}

\newcommand{\re}{\mathbb R}

\newcommand{\<}{\left<}

\newcommand{\ds}{\displaystyle}
\renewcommand{\>}{\right>}
\newcommand{\flecha}{\longrightarrow}

\newcommand{\gorro}{\widetilde}

\newcommand{\mc}[1]{\mathcal{#1}}

\input xypic.sty

\newcommand{\kh}{\ds \frac{K}{H^n}}

\def\det{{\rm det}}

\def\vle{{ Vol}}

\def\tr{{\rm tr}}

\def\fle{\rightarrow}
\def\parcial#1#2{\fracc{\partial #1}{\partial#2}}
\def\deri#1#2{\fracc{d #1}{d#2}}

\def\({\left (}
\def \){\right)}

\newcommand{\eps}{\ensuremath{\varepsilon}}

\newcommand{\usb}[2]{\!\!\!\! \underset{\eqref{#1}}{#2}}

\newtheorem{teor}{\hspace{12pt} Theorem}
\newtheorem{defi}[teor]{\hspace{12pt} Definition}

\newtheorem{prop}[teor]{\hspace{12pt} Proposition}
\newtheorem{lema}[teor]{\hspace{12pt} Lemma}

\newtheorem{notad}[teor]{\hspace{12pt} Remark}
\newtheorem{notap}[teor]{\hspace{12pt} Remark}
\newtheorem{coro}[teor]{\hspace{12pt} Corollary}

\numberwithin{teor}{section}

\newcommand{\be}{\begin{enumerate}}
\newcommand{\ee}{\end{enumerate}}
\newcommand{\bi}{\begin{itemize}}
\newcommand{\ei}{\end{itemize}}
\newcommand{\bd}{\begin{description}}
\newcommand{\ed}{\end{description}}
\newcommand{\bec}{\begin{equation}}
\newcommand{\eec}{\end{equation}}
\newcommand{\ba}{\begin{array}}
\newcommand{\ea}{\end{array}}
\newcommand{\bt}{\begin{teor}}
\newcommand{\et}{\end{teor}}
\newcommand{\bdem}{\begin{proof}}
\newcommand{\edem}{\end{proof}}
\newcommand{\bl}{\begin{lema}}
\newcommand{\el}{\end{lema}}
\newcommand{\bnp}{\begin{notap}}
\newcommand{\enp}{\end{notap}}
\newcommand{\bde}{\begin{defi}}
\newcommand{\ede}{\end{defi}}
\newcommand{\bnod}{\begin{notad}}
\newcommand{\enod}{\end{notad}}
\newcommand{\bp}{\begin{prop}}
\newcommand{\ep}{\end{prop}}
\newcommand{\bco}{\begin{coro}}
\newcommand{\eco}{\end{coro}}

\newcommand{\nn}{\nonumber}
\newcommand{\lb}{\label}

\newcommand{\ccdot}{\, \cdot \,}
\newcommand{\us}{\underset}
\newcommand{\fp}{\frak p}

\usepackage{bm}

\usepackage[draft]{optional}

\input xypic.sty

\renewcommand{\(}{\left(}

\renewcommand{\>}{\right>}
\renewcommand{\)}{\right)}
\def\bal{\begin{align}}
\def\eal{\end{align}}

\numberwithin{equation}{section}
\def\be{\begin{equation}}
\def\ee{\end{equation}}

\def\a{\alpha}

\def\tr{{\rm tr}}
\def\det{{\rm det}}
\def\dist{{\rm dist}}

\def\vle{{\rm vol}}
\def\parcial#1#2{\frac{\partial #1}{\partial#2}}

\def\deri#1#2{\frac{d #1}{d#2}}

\def\flecha{\longrightarrow}
\def\fle{\rightarrow}

\def\ds{\displaystyle}

\begin{document}

\markright{}


\title{Volume\,-preserving flow by powers of the $m^{\text{th}}$ mean curvature}

\pagestyle{myheadings}

\title{Volume\,-preserving flow \\ by powers of the $m^{\text{th}}$ mean curvature}

\author{ Esther Cabezas-Rivas \medskip \\
\footnotesize{Mathematics Institute, University of Warwick} \\
\footnotesize{Coventry CV4 7AL, United Kingdom } \\
\footnotesize{e-mail: E.Cabezas@warwick.ac.uk} \bigskip \\
Carlo Sinestrari \medskip \\
\footnotesize{Universit\`a di Roma ``Tor Vergata''} \\
\footnotesize{Via della Ricerca Scientifica 1, 00133 Roma, Italy} \\
\footnotesize{e-mail: sinestra@mat.uniroma2.it}}

\date{}

\maketitle

\begin{abstract}
We consider the evolution of a closed convex hypersurface under a volume preserving curvature flow. The speed is given by a power of the $m^{\text{th}}$ mean curvature plus a volume preserving term, including the case of powers of the mean curvature or of the Gauss curvature. We prove that if the initial hypersurface satisfies a suitable pinching condition, the solution exists for all times and converges to a round sphere.
\end{abstract}

\section{Introduction and main results}
\label{intro}

Let $M$ be a closed oriented $n$-dimensional manifold, with $n \geq 2$, and let $X_0 : M \flecha \re^{n +1}$ be a smooth immersion of $M$ into the euclidean space. 
We consider a family of immersions $X: M \times [0,T) \flecha \re^{n +1}$  which satisfies  the partial differential equation
\bec \lb{vpgcf}
\parcial{X}{t}(\ccdot, t) = \big(h(t) - \sigma(\cdot, t)\big) N_t, \qquad t \in [0,T)
\eec
with initial value
\bec \lb{inval}
X(\ccdot, 0) = X_0(\ccdot).
\eec
Here 
\bi
\item $N_t$ is the unit normal vector field along the immersion which induces on $M$ the given positive orientation. When $X(M)$ encloses a compact domain $\Omega$, the orientation is chosen so that $N$ points outward; 

\item $\sigma$ is a symmetric function of the principal curvatures of $M_t:=X(M,t)$;

\item $h(t)$ stands for  the averaged $\sigma$:  
\bec \ds 
h(t) = \ds \frac{1}{|M_t|} \int_M \sigma\, d \mu_t \lb{def_ht}, \eec
where $|M_t| := \int_M d \mu_t$ gives the area ($n$-dimensional volume) of $M_t$. With this definition, the flow \eqref{vpgcf} preserves the volume of the domain $\Omega_t$ enclosed by $M_t$ (when such an $\Omega_t$  exists).
\ei

In order to specify the class of speeds $\sigma$ we are going to consider, let us introduce some notation. We set $X_t:=X(\cdot,t)$ and we denote by $M_t$  both the immersion $X_t:M\flecha \re^{n + 1}$ and the image $X_t(M)$, as well as the Riemannian manifold $(M,g_t)$ with the metric $g_t$ induced by the immersion. We call $k_{1} \leq k_{2} \leq \ldots  \leq k_{n}$ the principal curvatures of $M_t$. We use the letters $H$ and $K$ for the mean curvature and Gauss curvature respectively, i.e. $H=k_1+\cdots+k_n$ and $K=k_1 \cdots k_n$. In addition, for any integer $m \in \{1,\dots,n\}$, we denote by $H_m$ the $m^{\text{th}}$ mean curvature, defined as
\bec \lb{def_mmc}
H_m = \frac{m! \, (n-m)!}{n!}\sum_{1 \leq i_1 < \ldots < i_m \leq n} k_{i_1} k_{i_2} \cdots k_{i_m}.
\eec 
Observe that $H_1=H/n$ and $H_n=K$; in addition, $H_2$ coincides, up to a constant factor, with the scalar curvature. Thus, the $m^{\text{th}}$ mean curvatures can be regarded as generalizations of these quantities. Various problems involving these curvatures have been considered in the literature on geometric analysis, such as finding hypersurfaces with prescribed $H_m$ curvature, see e.g. \cite{CNS4,GG}.

In this paper we consider the flow \eqref{vpgcf} with the speed $\sigma$ given by a power of an $m^{\text{th}}$ mean curvature, namely
\bec \lb{def_sig}
\sigma (k_1,\dots,k_n)= H_m(k_1,\dots,k_n)^\beta
\eec 
for some $\beta > 1/m$. In this way $\sigma$ is a homogeneous function of the curvatures with a degree $m \beta>1$. We call the flow \eqref{vpgcf} with this choice of speed {\bf volume-preserving flow by powers of the $m^{\text{th}}$ mean curvature}. Our analysis is focused on the behavior of convex hypersurfaces. The main results of our paper are summarized in the statement below.

\bt \lb{t_main}Given $m \in \{1,\dots,n\}$ and $\beta > 1/m$ there exists a constant $C_\fp=C_\fp(n,m,\beta) \in \,(0,1/n^n)$ such that, if the initial immersion $X_0$ satisfies at every point
\bec \lb{in_pinch}
K > C_\fp H^n >0,
\eec 
then problem \eqref{vpgcf}--\eqref{inval}, with $\sigma$ given by \eqref{def_sig}, has a unique solution, which satisfies the following properties:
\bi
\item[\rm{(a)}] Inequality \eqref{in_pinch} holds everywhere on $M_t$ for all $t>0$ such that the flow exists.

\item[\rm{(b)}] $M_t$ exists for all  $t \in [0, \infty)$.

\item[\rm{(c)}] The $M_t$'s converge, exponentially in the $C^\infty$-topology, to a round sphere enclosing the same volume as $M_0$.
\ei
\et

It is easy to check (see Section 4) that if $K > CH^n>0$ on a closed hypersurface then the principal curvatures are positive everywhere and satisfy $k_1 \geq \eps k_n$ for a suitable $\eps(C)>0$ which increases with $C$. Thus, condition \eqref{in_pinch} implies in particular the convexity of $M_t$, but it is a stronger assumption; it can be regarded as a pinching condition on the curvatures. Whereas we are not able to show that convexity is invariant under the flow \eqref{vpgcf}, we can prove that condition \eqref{in_pinch} is preserved by the flow for a suitable $C_\fp$. A similar pinching condition has been considered by Chow \cite{Ch_K} and Schulze \cite{Sch2} in the analysis of other flows.

There exists a wide literature about the behavior of convex hypersurfaces under geometric flows of the form \eqref{vpgcf} without the volume-preserving term $h(t)$. A classical result by Huisken \cite{Hu84} states that any closed convex hypersurface evolving by mean curvature flow shrinks to a point in finite time, and it converges to a sphere after an appropriate rescaling. This result was soon followed by similar ones due to K. Tso and B. Chow, in the case where $\sigma$ is a power of the Gauss curvature \cite{Tso, Ch_K} or of the scalar curvature \cite{Ch3}. After this, many other cases have been investigated where $\sigma$ is a homogeneous symmetric function of the principal curvatures.

Let us remark that an important role in the analysis of the case $h(t) \equiv 0$ is played by the degree of homogeneity of $\sigma$. In fact, the result of \cite{Hu84} has been generalized to a large class of $\sigma$ homogeneous of degree one (see in particular \cite{An94}). When the degree is greater than one, the analysis is less complete and often restricted to dimension two, see \cite{roll,A04,OS_1}; the results in general dimension concern specific choices of the speed  \cite{Ch_K, Sch2, AS} and all of them require a pinching condition on the initial hypersurface. The case where the degree is less than one is even more difficult. In some case it is known that convex hypersurfaces shrink to a point  \cite{Agcf,Sch1}, but some counterexamples show that in general the profile may not become spherical after rescaling.

It is natural to investigate how this behavior changes when we consider the volume-preserving version of these flows, that is, when we add the $h(t)$ term in equation \eqref{vpgcf}. Since in this case surfaces cannot shrink to a point, one expects that the flow starting from convex data exists for all times and converges smoothly to a round sphere. In the case of the volume-preserving mean curvature flow, this property was proved, again by Huisken, in \cite{Hu87}. By a different approach, Escher and Simonett \cite{EsSi98} proved the same property when the initial hypersurface is a small perturbation of the sphere (possibly non-convex). Andrews \cite{An01} studied the anisotropic version of this flow, while McCoy \cite{McC1,MC04} studied mean curvature flows which preserve other kinds of volume. Recently, the first author and Miquel \cite{CaMi1} considered volume-preserving mean curvature flow in hyperbolic space and proved convergence to geodesic spheres.

On the other hand, speeds different from the mean curvature have been treated only by McCoy \cite{Mc3} who proved the convergence to a sphere  for a large class of functions $\sigma$ homogeneous of degree one (including the case $\sigma=H_m^\beta$ for $\beta=1/m$). In this paper we consider instead speeds with a larger degree of homogeneity. The class given by  \eqref{def_sig} is a natural choice, because on one hand it includes the most significant examples and, on the other hand, it allows to use some properties of the elementary symmetric polynomials (parts (b) and (e) of Lemma \ref{lpolS}) which are essential in the analysis of the PDEs associated to the flow.

As in most of the literature quoted above, a fundamental step in our procedure consists of finding a pinching condition which is invariant under the flow (part (a) of Theorem \ref{t_main}). For the class of speeds we consider, we cannot follow the procedure of \cite{Mc3}; instead, we apply the maximum principle to the quotient $K/H^n$, after a careful estimation of the gradient terms in the equation. Once the preservation of the pinching condition has been established, we can prove by similar arguments as in \cite{Mc3} that the curvature remains uniformly bounded as long as the flow exists.

To prove long time existence of the flow (part (b) of Theorem \ref{t_main}), we then need uniform bounds on all derivatives of the curvature. This is equivalent to obtaining higher order derivative estimates for a suitable fully nonlinear parabolic equation. Since the homogeneity of the speed is greater than one, the operator is not concave in the second derivatives and the usual theory by Krylov and Safonov \cite{K,Lieb} cannot be applied. We use instead a technique introduced by Andrews in \cite{A5} (see also \cite{Tsai}) based on a regularity result by Caffarelli \cite{Caf1} for elliptic equations. 

Another complication due to the higher degree of homogeneity, similar to the one which was pointed out in \cite{Sch2}, occurs in the analysis of the asymptotic behavior of our flow (part (c) of Theorem \ref{t_main}). Since the coefficients of the second order operator depend on the curvature, the equation may a priori degenerate when time goes to infinity because we lack a positive lower bound on the curvature
which is uniform in time. To deal with this problem, we first estimate from below the possible decay rate of the curvature. Then we rewrite the equation satisfied by the speed as a suitable porous medium equation. Finally, as in \cite{Sch2}, we prove convergence to a sphere after applying a regularity result on degenerate parabolic equations due to Di Benedetto and Friedman \cite{DiB_F}.

Some of the results of this paper, corresponding to the case $\sigma=H_n^\beta$ (i.e. a power of the Gauss curvature), were already obtained in \cite{tesis}.

\section{Preliminaries}

Throughout the paper, we shall use $\bar g = \<\ccdot, \ccdot\>$ and $\overline \nabla$ to represent the metric and covariant derivative, respectively, of $\re^{n + 1}$. Then each immersion $X_t$ of $M$ into $\re^{n + 1}$ induces the following symmetric 2-covariant tensor fields:

$\bullet$ a {\it metric} $g_t$ defined by $g = X^\ast \bar g$, and

$\bullet$ a {\it second fundamental form} $\a_t$ given by $\a(\ccdot, \ccdot) = \<\overline \nabla_{X_\ast \ccdot} N, X _\ast \ccdot\>.$

From here, we can also introduce the {\it shape operator} $A$ (or {\it Weingarten map}) of $X$ as  $\a(\ccdot, \ccdot) = g(A\ccdot, \ccdot)$; the eigenvalues $k_1\leq \ldots \leq k_n$ of $A$ are called the principal curvatures. We say that $M_t$ is {\it convex} if $k_1 \geq 0$ everywhere and that it is {\it uniformly convex} if $k_1>0$ everywhere. The {\it mean curvature} is given by $H = \tr_g \a = \sum_i k_i$ and the {\it Gauss curvature} by $K = \frac{\det \a}{\det g} = \prod_i k_i$. More generally, we call m$^{\text{th}}$ mean curvature of a hypersurface the function $H_m$ defined in \eqref{def_mmc}, that is, the m$^{\text{th}}$ elementary symmetric polynomial of the principal curvatures, up to a normalizing factor. Since $H_m$ is homogeneous of degree $m$, the speed $\sigma$ is homogeneous of degree $m \beta$ in the curvatures $k_i$.

We shall often use the symbol $\frak k$ to denote the vector $(k_1,k_2,\dots,k_n)$ whose entries are the principal curvatures or, depending on the context, a generic element of $\re^n$. We denote by $\Gamma_+ \subset \re^n$ the positive cone, i.e.
$$\Gamma_+=\{ {\frak k}=(k_1,\dots,k_n) ~:~ k_i>0 \mbox{ for all }i\}.$$
While the polynomials $H_m$ are defined for any ${\frak k} \in \re^n$, $\sigma$ is well defined and smooth for a general $\beta$ only on the cone where $H_m$ is positive; such a cone clearly includes $\Gamma_+$ for any $m$.

Observe that $H,K,H_m,\sigma$, may be regarded as functions of $\frak k$, or as functions of $A$, or as functions of $\alpha,g$, or also as functions of space and time on $M_t$. For the sake of simplicity, we denote these functions by the same letters in all cases, since the meaning should be clear from the context. We use the notation
\bec\lb{e.sigma}
\dot \sigma^i :=\frac{\partial \sigma}{\partial k_i}, \quad \text{and} \quad \dot \sigma^{ij} :=\frac{\partial \sigma}{\partial \alpha_{ij}},
\eec
and also
$$
\tr(\dot \sigma) :=  \sum_{i=1}^{n} \frac{\partial \sigma}{\partial k_i} =\frac{\partial \sigma}{\partial \alpha_{ij}} g^{ij}.
$$
In addition, if $B,\bar B$ are matrices, we write
$$\dot \sigma(B) := \parcial{\sigma}{\a_{ij}} B_{ij} \quad \text{and} \quad \ddot \sigma(B, \bar B) = \frac{\partial^2 \sigma}{\partial \a_{ij} \, \partial \a_{ls}} B_{ij} \bar B_{ls}.$$
Hereafter, we use the Einstein convention of sum over repeated indices.

\subsection{Basic properties of the m$^{\text{\bf th}}$ mean curvatures}

The following lemma gathers some properties of $H_m$ which will be used repeatedly throughout the paper.

\bl \lb{lpolS} For any $m \in \{1, \ldots, n\}$, the function $H_m$ satisfies: 
\bi 
\item[{\rm (a)}]   $\ds \parcial{H_m}{k_i} ({\frak k})>0$ for all $i \in \{1, \ldots, n\}$ and ${\frak k} \in \Gamma_+$.

\item[{\rm (b)}] $H_m^{1/m}$ is concave in $\Gamma_+$.

\item[{\rm (c)}]  $\tr(\dot \sigma) \geq m \, \beta \, \sigma^{1 -  \frac{1}{m\beta}}$.

 \item[{\rm (d)}]  $H_m^{1/m} \leq \ds\frac{H}{n}$; equivalently, $\sigma \leq \(\ds\frac{H}{n}\)^{m \beta}$.

\item[{\rm (e)}]  $H_m$, as a function of $\a_{ij}$, is also a homogeneous polynomial of degree $m$; in addition, as a function on $M$, it satisfies $\nabla_j\(\ds\parcial{H_m}{\a_{ij}}\) = 0$ for any $i \in \{1, \ldots, n\}$, where $\nabla$ is the covariant derivative on $M$.

 \ei
\el

\bdem
Property (a) is a direct consequence of the definition of $H_m$. Part (b) is a well known property (see e.g. \cite[Theorem 15.16]{Lieb}). Statements (c) and (d) follow from Lemma 4.2 and estimate (3.16) in \cite{Ur1} respectively. Property (e) is also well known (see e.g. Proposition 2.1 (a)-(b) and Proposition 3.3 in \cite{Rei}). \edem

\subsection{Some estimates for convex hypersurfaces}

We state here three auxiliary results on convex hypersurfaces. Henceforth we shall  denote by $\rho$ the {\it inradius} of $M = \partial \Omega$ and by $D$ the {\it outer radius}, which are respectively the radius of the biggest ball enclosed by $M$ and of the smallest ball which encloses $M$.

The first result was observed by Andrews in \cite{An94}, \!\! Lemma 5.4 and Theorem 5.1; it shows that a pinching inequality on the curvatures implies a bound on the ratio between outer radius and inradius.

\bl\lb{lem_An}
Let $X: M^n \fle \re^{n +1}$ be a smooth, uniformly convex immersion of the compact manifold $M^n$. If  
$$k_n (x) \leq B_1 k_1(x)$$
for every $x \in M^n$ and for some constant $B_1 < \infty$, then we have 
$$D \leq B_2 \, \rho$$
for another constant $B_2$ depending on $n$ and $B_1$.\el

The next result is a particular case of the  {\em Alexandrov-Fenchel inequality} for convex sets (see e.g. \cite{McC1}, pp. 28--29, and the references therein).

\bl \lb{Alex_Fen}
There exists a constant $C_n > 0$ such that, for any bounded smooth convex region $\Omega \subset \re^{n + 1}$, we have 
$$\(\int_{\partial \Omega} H \, d \mu\)^{n + 1} \geq C_n \, \vle(\Omega)^{n - 1}.$$
\el

Finally, we recall an algebraic property proved by Schulze in \cite[Lemma 2.5]{Sch2}.

\bl \lb{2.5} For any $\eps>0$ there exists $\delta=\delta(\eps,n)>0$ with the following property: if we have $k_1 \geq \eps H > 0$ at some point of an $n$-dimensional hypersurface, then at the same point we also have
\bec \lb{Sch_2.5}
\frac{n |A|^2 - H^2}{H^2} \geq \delta \(\frac1{n^n} - \kh\).
\eec 
\el

\subsection{H\"older estimates for nonlinear PDEs}  \lb{par_teor}

In our analysis we need some a priori estimates on the H\"older norms of the solutions to elliptic and parabolic partial differential equations in euclidean spaces. We recall that, in the case of a function depending on space and time, there is a suitable definition of H\"older norm which is adapted to the purposes of parabolic equations (see e.g. \cite{Lieb}). In addition to the standard Schauder estimates for linear equations, we use in the paper some more recent results which are collected here.  The estimates below hold for suitable classes of weak solutions; for the sake of simplicity, we state them in the case of a smooth classical solution, which is enough for our purposes.

Given $r>0$, we denote by $B_r$ the ball of radius $r>0$ in $\re^n$ centered at the origin. First we recall a well known result due to Krylov and Safonov, which applies to linear parabolic equations of the form
\bec \lb{ec_par}
\(a^{ij}(x, t) D_{i} D_j  + b^i(x, t) D_i  + c(x, t) - \parcial{}{t}\) u = f
\eec
in $B_r \times [0,T]$, for some $T>0$.
We assume that $a^{ij} = a^{ji}$ and that $a^{ij}$ is elliptic; that is, there exist two constants $\lambda, \, \Lambda > 0$ such that
\bec \lb{SP}
\lambda |v|^2 \leq a^{ij}(x,t)v_iv_j \leq \Lambda |v|^2
\eec
for all $v \in \re^n$ and all $(x, t) \in B_r \times [0,T]$. Then the following estimate holds \cite[Theorem 4.3]{KS}:

\bt \lb{KS_t}
Let $u \in C^2(B_r \times [0,T])$ be a solution of \eqref{ec_par}, where the coefficients are measurable, satisfy \eqref{SP} and 
$$|b^i|, \, |c|  \leq K_1 \qquad \text{for all } i = 1, \ldots, n,$$
for some $K_1>0$. Then, for any  $0<r'<r$ and any $0<\delta <T$ we have
$$\|u\|_{C^{\a}(B_{r '} \times [\delta, T])} \leq C\(\|u\|_{C(B_r \times [0,T])} + \|f\|_{L_\infty(B_r \times [0,T])}\)$$
for some constants $C > 0$ and $\a \in (0,1)$ depending on $n$, $\lambda$, $\Lambda$, $K_1$, $r$, $r'$ and $\delta$.
\et

Next we quote a result for fully nonlinear elliptic equations, which is due to Caffarelli. We consider the equation
\bec \lb{FNL}
F(D^2u(x),x)=f(x), \qquad x \in B_r.
\eec
Here $F:{\cal S} \times B_r \to \re$, where $\cal S$ is the set of the symmetric $n \times n$ matrices. The nonlinear operator $F$ is called elliptic if there exist $\Lambda \geq \lambda >0$ such that
\bec \lb{FNL-ell}
\lambda ||B|| \leq F(A+B,x) - F(A,x) \leq \Lambda ||B||
\eec
for any $x \in B_r$ and any pair $A,B \in \cal S$ such that $B$ is nonnegative definite.

\bt  \lb{Caf-pert}
Let $u \in C^2(B_r)$ be a solution of \eqref{FNL}, where $F$ is continuous and satisfies \eqref{FNL-ell}. Suppose in addition that $F$ is concave with respect to $D^2u$ for any $x \in B_r$. Then there exists $\bar \alpha \in \,(0,1)$ with the following property: if, for some $K_2>0$ and $\alpha \in \,(0,\bar \alpha)\,$, we have that $f \in C^\alpha(\Omega)$ and that
$$
F(A,x) - F(A,y) \leq K_2 |x-y|^{\alpha} (||A||+1), \qquad x,y \in B_r, \  A \in \cal S
$$
then, for any $0<r'<r$, we have the estimate
$$\|u\|_{C^{2+\a}(B_{r'})} \leq C(||u||_{C(B_r)}+||f||_{C^\alpha(B_r)}+1)$$
where $C > 0$ only depends on $n$, $\lambda$, $\Lambda$, $K_2$, $r$ and $r'$.
\et

The above result follows from Theorem 3 in \cite{Caf1} (see also Theorem 8.1 in \cite{CC} and the remarks thereafter). It generalizes, by a perturbation method, a previous estimate, due to Evans and Krylov, about equations with concave dependence on the hessian. In contrast with Evans-Krylov result (see e.g. inequality (17.42) in \cite{GT}), Theorem \ref{Caf-pert} gives an estimate in terms of the $C^\alpha$-norm of $f$ rather than the $C^2$-norm, and this is essential for our purposes.

Finally, we recall an interior H\"older estimate, due to Di Benedetto and Friedman \cite[Theorem 1.3]{DiB_F}, for solutions of the degenerate parabolic equation 
\bec \lb{DiBF_ec}
\parcial{v}{t} - D_i\(a^{ij}(x, t, D v) D_j v^d\) = f(x,t,v, D v),
\eec
being $d > 1$.

\bt \lb{tBF}
Let $v \in C^2(B_r \times [0,T])$ be a nonnegative solution of \eqref{DiBF_ec}, where $a^{ij}$ satisfy \eqref{SP}. Let $c_1, c_2, N>0$ be such that 
$$|f(x,t,v, D v)|\leq c_1 |D v^d| + c_2,$$
and
 $$\sup_{0 < t < T} ||v(\, \cdot \,, t)||_{L^2(B_r)}^2 + \|D v^d\|_{L^2(B_r \times [0,T])}^2 \leq N.$$ Then for any $0<\delta<T$ and $0<r'<r$, we have
$$\|v\|_{C^{\a}\(B_{r'} \times [\delta, T]\)}  \leq C,$$
for suitable $C>0, \alpha \in (0, 1)$ depending only on $n, N, \lambda, \Lambda, \delta, c_1, c_2, r$ and $r'$.
\et

\section{Short time existence and evolution equations}

As far as short time existence is concerned, the properties of the volume-preserving flows do not differ substantially from the ones of the standard flows without the volume-preserving term, as the result below shows. 

\bt \lb{ste_kb} Let $X_0: M \fle \re^{n +1}$ be a smooth, closed and uniformly convex hypersurface. Then there exists  a unique smooth solution $X(\ccdot, t)$ of  problem \eqref{vpgcf}--\eqref{inval}, defined on some time interval $[0, T)\,$, with $T > 0$.
\et

\bdem
We first note that the above result holds for the flow \eqref{vpgcf} without the volume-preserving term $h(t)$. In fact, it is well known (see e.g. \cite[\S 3]{HuiPol})  that a flow of the form
\bec \lb{gen-flow}
\parcial{X}{t}(\ccdot, t) =  - f(\cdot, t) N_t,
\eec
where $f$ is any symmetric function of the curvatures $k_i$, is parabolic on a given hypersurface if the condition $\parcial{f}{k_i} > 0$, with $i=1,\dots,n$, holds at every point. Then, given any initial immersion $X_0$ satisfying the parabolicity assumption, standard techniques ensure the local existence and uniqueness of a solution to \eqref{gen-flow} with initial value $X_0$. In our case we have $f=\sigma=H_m^\beta$ and the condition reads
\bec \lb{deri_sig1}
\parcial{\sigma}{k_i} = \beta H_m^{\beta -1} \parcial{H_m}{k_i} >0,
\eec
which is satisfied on any uniformly convex hypersurface, by Lemma \ref{lpolS} (a).

It is not difficult to extend this result to the volume-preserving flow \eqref{vpgcf}. As it is pointed out in \cite{Hu87},  the $h(t)$ term does not influence the parabolicity of the equation and so one can repeat the proof of the previous case with minor modifications. A more detailed argument is given by McCoy in \cite[\S 7]{Mc3}; although the assumptions on the speed in that paper are different, the proof applies to our case as well. \edem

\bnp
A closed hypersurface satisfying the pinching assumption \eqref{in_pinch} for some $C_\fp>0$ is uniformly convex. In fact, any closed hypersurface has at least one point where all curvatures are positive. Since $K>0$ by \eqref{in_pinch}, the curvatures cannot vanish and therefore are positive everywhere. This shows that the above existence result applies to the class of initial data considered in Theorem \ref{t_main}. 
\enp

The rest of this section is devoted to collect some basic evolution formulae under the flow \eqref{vpgcf}. The following lemma, as Lemma \ref{ev_eq2} below, can be obtained by computations similar to those in Section 3 of \cite{Hu84} (see also \cite{Ch_K,Mc3,tesis}). 

\bl If $M_t$ is a solution of \eqref{vpgcf}, the following evolution equations hold:
$$\ba{ll}
 \text{\rm (a)}  \ \ds \parcial{g}{t} = 2 (h -  \sigma)  \a \qquad  
 &  \text{\rm  (b)}  \ \ds \parcial{g^{ij}}{t} = - 2 (h -  \sigma) g^{il} g^{js} \a_{ls} \medskip \\
\text{\rm (c)}  \ \ds \parcial{N}{t} =  \nabla \sigma & \medskip \\
 \text{\rm (d)}  \ \ds \parcial{}{t} d\mu_t =  (h - \sigma) \ H \ d\mu_t \qquad & \text{\rm (e)} \ \ds \deri{}{t} |M_t| =  \int_M (h - \sigma) \ H \ d\mu_t,
\ea$$
where $\nabla$ denotes the covariant derivative induced by $g$.
\el

Notice that the flow is defined so that we get the volume-preserving property:
\bl 
The enclosed volume satisfies $V:= \vle(\Omega_0) = \vle(\Omega_t)$
for all $t$ such that the solution of \eqref{vpgcf} is well-defined, i.e. 
$$\deri{}{t} \vle(\Omega_t) = 0, \qquad
\text{ where }\partial \Omega_t = M_t.$$
\el

\bdem
As in \cite{tesis}, we use $\deri{}{t} \vle(\Omega_t) =  \int_M \<\parcial{X}{t}, N\> \, d\mu_{t} = \int_M (h - \sigma) \, d\mu_{t} = 0$.
\edem

Hereafter, given any $(2,0)$-tensor $a$, we shall denote
$$\Delta_a:= a^{ij} \nabla_i \nabla_j \qquad \text{ and } \qquad |B|^2_a := a^{ij} B_i B_j$$
(also $\Delta = \Delta_{g^{-1}}$ and $|\ccdot|^2 = |\ccdot|^2_{g^{-1}}$). A tensor which occurs repeatedly in the analysis of our flow is $\dot \sigma=(\dot \sigma^{ij})$ introduced in \eqref{e.sigma}. The following equations also follow from a direct computation.

\bl \lb{ev_eq2}
Under \eqref{vpgcf}, we have the evolution formulae listed below.

\medskip
{\rm (a)} 
\vspace*{-0.6cm}
\bec \lb{evF}
\hspace*{-5.8cm} \parcial{\sigma}{t} = \Delta_{\dot \sigma} \sigma + (\sigma - h) \,\tr_{\dot \sigma} (\a A).
\eec

\medskip

{\rm (b)} 
\vspace*{-1cm}
 \begin{align*}
\hspace*{-0.45cm}\parcial{\a}{t} &= \Delta_{\dot \sigma} \a  +  \ddot{\sigma}(\nabla_{_{\! \ds \cdot}} \a, \nabla_{_{\! \ds \cdot}} \a)   + \tr_{\dot \sigma}(\a A) \, \a + \big(h - (m \beta + 1) \sigma\big) \a A.  
\end{align*}

\medskip
{\rm (c)}
\vspace*{-1cm} 
\begin{align*} 
\qquad \parcial{H}{t} & = \Delta_{\dot \sigma} H + \tr_{g^{-1}} \big[\ddot{\sigma}(\nabla_{_{\! \ds \cdot}} \a, \nabla_{_{\! \ds \cdot}} \a)\big]  - \big(h + (m \beta -1) \sigma\big) |A|^2 
 + H \tr_{\dot{\sigma}}(\a A).
\end{align*}

{\rm (d)} \vspace*{-1cm} 
\begin{align*} 
\parcial{K}{t} & = \Delta_{\dot \sigma} K   - \frac{|\nabla K|_{\dot \sigma}^2}{K} + K \(\tr_b \big[\ddot{\sigma}(\nabla_{_{\! \ds \cdot}} \a, \nabla_{_{\! \ds \cdot}} \a)\big]  - \tr_{\dot \sigma}\(\nabla_{_{\! \ds \cdot}}  b^{ij} \nabla_{_{\! \ds \cdot}}  \a_{ij}\)\)\nn
\\ & \quad  - \big(h + (m \beta -1) \sigma\big) H K
 + n K \tr_{\dot{\sigma}}(\a A),
\end{align*}
where $b := \a^{-1}$.

{\rm (e)} For the position vector field $\frak X(\ccdot, t) - \bar x$  (with origin $\bar x$) on $M_t$: 
\bec \lb{fal_sopev}
\parcial{}{t} \<\frak X, N\> = \Delta_{\dot \sigma} \<\frak X, N\> + \tr_{\dot \sigma} (\a A) \<\frak X, N\> + \big(h - (m \beta + 1) \sigma\big).
\eec
\el

Next, with the purpose of writing the evolution equation for $H_m$, we set
\bec \lb{def_c}
c^{ij} = \frac{\partial H_m}{\partial \a_{ij}}. 
\eec
Recall (cf.  part (e) of Lemma \ref{lpolS}) that the symmetric tensor $c$ is divergence-free. With this notation, we obtain
\bl
If $M_t$ is a hypersurface in $\re^{n + 1}$ evolving under \eqref{vpgcf}, the $m^{\text{th}}$ mean curvature $H_m$ and its $\beta^{\text{th}}$ power $\sigma$ satisfy the following evolution equations:
$$\parcial{H_m}{t} = \beta H_m^{\beta -1}\(\Delta_c H_m + (\beta -1) \frac{|\nabla H_m|_c^2}{H_m}\) + (\sigma - h) \tr_{c} (\a A) 
$$ 
and
\bec
\parcial{\sigma}{t} = \beta H_m^{\beta -1} \Delta_c \sigma + \beta \frac{\sigma - h}{H_m} \, \tr_c (\a A) \,\sigma. \lb{ev_sig2}
\eec
\el

\bdem
From \eqref{def_c}, we can write $\dot \sigma^{ij} = \beta H_m^{\beta - 1} c^{ij}$; thus the evolution equation for $H_m$ becomes
\begin{align}
\parcial{H_m}{t} = \frac1{\beta H_m^{\beta - 1}} \parcial{\sigma}{t} \us{\eqref{evF}}{=} \Delta_c \sigma +  (\sigma - h) \tr_{c}(\a A). \lb{Hm_pm}
\end{align}
So the first formula in the statement follows by applying the product rule for $\Delta_c$. Arguing in a similar way, we obtain \eqref{ev_sig2}.
\edem

In the next statement we use the notation $|T|^2_{a,b}:= a^{ij} b^{rs} b^{lp} T_{isl} T_{jrp}$ for any 3-tensor $T$.
\bl The quantity $q:= K/H^n$ evolves under \eqref{vpgcf} satisfying
\begin{align}
\parcial{q}{t}  & = \Delta_{\dot \sigma} q + \frac{(n+1)}{n H^n} \<\nabla H^n, \nabla q\>_{\dot \sigma} - \frac{(n-1)}{n K} \< \nabla K, \nabla q\>_{\dot \sigma}  - \frac{H^n}{n K} \left|\nabla q\right|^2_{\dot \sigma} \nn
\\ & \quad+ \frac{q}{H^{2}} \left|H \nabla \a - \nabla H \, \a\right|^2_{\dot \sigma,b} + q\,\tr_{b - \frac{n}{H} g^{-1}} \big[\ddot{\sigma}^\flat(\nabla_{_{\! \ds \cdot}} \a, \nabla_{_{\! \ds \cdot}} \a)\big]  \nn
\\ & \quad + \(h + (m \beta -1) \sigma\) \frac{q}{H} \( |A|^2 n - H^2\).  \lb{pin_vpgcf}
\end{align}
\el

\bdem
The above formula follows from parts (c) and (d) of Lemma \ref{ev_eq2} by a straightforward computation similar to Lemma 2.2 in \cite{Sch2}. To express the gradient terms in the desired form, we use the identity
\begin{align*}
\tr_{\dot \sigma}\!\!\left[\frac{\nabla_{_{\! \ds \cdot}}  K \nabla_{_{\! \ds \cdot}}  K}{K^2}  + \nabla_{_{\! \ds \cdot}}  b^{kl} \nabla_{_{\! \ds \cdot}}  \a_{kl}\right] & =    \frac{(n -1) |\nabla K|_{\dot \sigma}^2}{n K^2} + \frac{H^{2n}}{n K^2} \left|\nabla \frac{K}{H^n}\right|^2_{\dot \sigma}
    - \frac{\left|H \nabla \a - \nabla H \, \a\right|^2_{\dot \sigma, b}}{H^2}
\end{align*}
which follows by taking traces with $\dot \sigma^{ij}$ in the formula at the bottom of  p. 121 in \cite{Ch_K}.
\edem

\section{The pinching estimate and its consequences}

In this section we prove a monotonicity property for the quotient $q=K/H^n$. Such a quotient, which was also considered in \cite{Ch_K, Sch2}, is a natural quantity to deal with in the study of our flow. Observe that, by the arithmetic-geometric mean inequality, we have $q \leq 1/n^n$, with equality only if $k_1=\dots=k_n$, i.e. at an umbilical point. Thus, the only hypersurfaces such that $q \equiv 1/n^n$ are the spheres. In addition, a lower bound of the form $K/H^n \geq C >0$ implies a pinching condition on the curvatures of the form $k_1 \geq \eps k_n$ which has various useful consequences for the analysis of our problem.

In order to apply the maximum principle to equation (\ref{pin_vpgcf}), we first derive some preliminary inequalities. The following result, pointed to us by G. Huisken, is a stronger version of Lemma 2.3 (ii) in \cite{Hu84}. 

\bl
If for some $\eps>0$ the inequality $\a \geq \eps Hg>0$ holds at a point of a hypersurface immersed in $\re^{n+1}$, then $\eps \leq 1/n$ and at the same point we also have
\bec
|H \nabla \a  - \a \nabla H|^2 \geq \frac{n-1}2\eps^2 H^2 |\nabla A|^2. \lb{leHuS} 
\eec
\el
\bdem
The assumption is equivalent to $k_1 \geq \eps H>0$. Using this, we first observe that
$$
H= k_1+\dots+k_n \geq n k_1 \geq n \eps H>0,
$$
which implies that $\eps \leq 1/n$. Again using that $k_1 \geq \eps H$, we also deduce
\begin{align}
A^j_i A^l_j \nabla^i H \nabla_l H & \leq k_n^2 |\nabla H|^2 \leq (|A|^2-(n-1) k_1^2) |\nabla H|^2
\nonumber \\
& \leq (|A|^2-(n-1) \eps^2 H^2 ) |\nabla H|^2  \leq  (1-(n-1) \eps^2) |A|^2 |\nabla H|^2. \lb{laux1}
\end{align}
Now we can write
\bec 
|H \nabla \a  - \a \nabla H|^2 = |\nabla A|^2 H^2 + |\nabla H|^2 |A|^2 - 2 H \<\a \nabla H, \nabla \a\> \lb{laux2}
\eec
whose last term, in local coordinates, is equal to $$-2 H \nabla_i \a_{jl} \nabla^i H \a^{jl} = - H \nabla_i \a_{jl}\(\nabla^i H \a^{jl} + \nabla^j H \a^{il}\),$$
by the Codazzi equations. Taking this into account and using the inequality $\<T, V\> \leq \frac{2 - \eps'}{2} |T|^2 + \frac1{2(2 - \eps')} |V|^2$, with $T=H \nabla_i \a_{jl}$, $V=
\nabla^i H \a^{jl} + \nabla^j H \a^{il}$ and $\eps' = (n-1)\eps^2$, we estimate
\begin{align*}
2 H \<\a \nabla H, \nabla \a\> & \leq \frac{2 - \eps'}{2} H^2 |\nabla A|^2 + \frac1{2(2 - \eps')} |\nabla^i H \a^{jl} + \nabla^j H \a^{il}|^2
\\& = \(1 - \frac{\eps'}{2}\) H^2 |\nabla A|^2 + \frac1{2 - \eps'}\(|\nabla H|^2 |A|^2 +  \nabla_i H \a_{jl} \nabla^j H \a^{il}\)\\ &\! \!\us{\eqref{laux1}}{\leq} \(1 - \frac{\eps'}{2}\) H^2 |\nabla A|^2 + |\nabla H|^2 |A|^2,
\end{align*}
which implies \eqref{leHuS} by means of \eqref{laux2}.
\edem

We also need the following elementary property.
\bl \lb{lem.const}
Given any $\eps \in (0,1/n)$, there exists $C=C(\eps,n) \in (0,1/n^n)$ such that, for any $\frak k=(k_1,\dots,k_n) \in \re^n$ with $k_i \geq 0$ for all $i=1,\dots, n$, we have
$$
K(\frak k) > C {H^n}(\frak k) \ \Longrightarrow \min_{1 \leq i \leq n} k_i > \eps H(\frak k),
$$
where $K(\frak k)=k_1\cdots k_n$ and $H(\frak k)=k_1+\dots+k_n$.
\el
\bdem
 Given any $\eps \in (0,1/n)$, we define
$$
{\cal A}_\eps= \{ \frak k =(k_1,\dots,k_n) ~:~  0 \leq   \min_{1 \leq i \leq n}k_i \leq \eps H(\frak k) \},
$$
$$
{\cal C}_\eps=\{ \frak k \in {\cal A}_\eps  ~:~ |\frak k |=1 \}.
$$ 
We have $H(\frak k)>0$ on any nonzero element of ${\cal A}_\eps$; hence, the quotient $K/H^n$ is defined everywhere on ${\cal C}_\eps$. Let us call $M_\eps$ the maximum of $K/H^n$ on ${\cal C}_\eps$, which exists because ${\cal C}_\eps$ is compact. Observe that $M_\eps < 1/n^n$. In fact, $K/H^n  \leq 1/n^n$, with equality if and only if $k_1=\dots=k_n$. Therefore, if $ \frak k$ is such that $K/H^n (\frak k)= 1/n^n$, then $\frak k$ satisfies $k_i=H(\frak k)/n$ for all $i$ and does not belong to ${\cal C}_\eps$ because we assume $\eps < 1/n$. 

By homogeneity, the inequality $K \leq M_\eps H^n$ is also satisfied by the elements of ${\cal A}_\eps$. Therefore, if $\frak k=(k_1,\dots,k_n)$ with $k_i \geq 0$ for all $i$ is such that $K > M_\eps H^n$, then $\frak k$ does not belong to ${\cal A}_\eps$.  The lemma follows by choosing $C=M_\eps$. 
\edem

We are now ready to prove a pinching estimate for our flow, which is one of the key steps in the proof of our main result.

\bt \lb{pinching}  
There exists a constant $C_\fp=C_\fp(n,m,\beta) \in (0,1/n^n)$ with the following property: if $X:M \times (0,T)\, \to \re^{n+1} $, with $t \in (0,T)$, is a smooth solution of \eqref{vpgcf}--\eqref{inval}, with $\sigma$ given by \eqref{def_sig} for some $\beta>1/m$, such that
\bi
\item the initial immersion $X_0$ satisfies \eqref{in_pinch} with the constant $C_\fp$,
\item the solution $M_t=X(M,t)$ satisfies $H>0$ for all times $t \in (0,T)$,
\ei
then the minimum of $K/H^n$ on $M_t$ is nondecreasing in time.
\et

\bnp We will see later in Corollary \ref{pos_H} that the above theorem is still valid without requiring that $H>0$ for $t \in (0,T)$.
\enp

\bdem The assumption that $H > 0$ ensures that the quotient $q(\ccdot, t)=K/H^n(\cdot, t)$ is well-defined for $t\in (0, T)$. Let us denote $\frak q(t) := \min_M q(\cdot, t)$.

Observe that it suffices to prove the theorem under the additional hypothesis that all principal curvatures are positive on $M_t$ for any $t \in (0,T)\,$. In fact, this holds for $t=0$ by \eqref{in_pinch}. If
there exists a first time $t_0>0$ at which $k_1=0$ at some point, we have $\frak q(t_0) = 0$. On the other hand, if the theorem holds in the convex case, $\frak q$ is nondecreasing in $(0, t_0)$; so it cannot decrease from $C_{\frak p}$ to zero. Thus, we can assume that the curvatures of $M_t$ are positive.

Now recall the well known fact (see e.g. \cite[\S 10.3]{ChII}) that 
$$D_+ \frak q(t) \geq \inf_{\mc M(t)} \parcial{q}{t}$$
where $\mc M (t) :=\{p \in M/  q(p, t) = \frak q(t)\}$ and $D_+$ denotes the lower right Dini derivative.
Since $\Delta_{\dot \sigma} q \leq 0$ and $\nabla q = 0$ on $\mc M(t)$, thanks to \eqref{pin_vpgcf} we get
\begin{align}
D_+ \frak q &\geq  \frak q\(\frac{1}{H^{2}} \left|H \nabla \a - \nabla H \, \a\right|^2_{\dot \sigma,b} + \,\tr_{b - \frac{n}{H} g^{-1}} \big[\ddot{\sigma}(\nabla_{_{\! \ds \cdot}} \a, \nabla_{_{\! \ds \cdot}} \a)\big]\)  \nn
\\ & \geq \frak q\(\frac{1}{H^{2}} \left|H \nabla \a - \nabla H \, \a\right|^2_{\dot \sigma,b} - \,\left|b - \frac{n}{H} g^{-1}\right| \big|\ddot{\sigma}(\nabla A, \nabla A)\big|\), \lb{auxP0}
\end{align}
where we have also used that the last term in
\eqref{pin_vpgcf} is nonnegative (by convexity and the elementary inequality $|A|^2 \geq H^2/n$). 

We want to show that the above expression is nonnegative provided the second fundamental form is suitably pinched. To this purpose we need to bound from below the positive term. Doing computations at a point where we choose an orthonormal basis which diagonalizes $\alpha$, we first deduce
\begin{align}
\left|H \nabla \a - \nabla H \, \a\right|^2_{\dot \sigma,b} & = \sum_{i,j,l} \dot \sigma^i \frac1{k_j}\frac1{k_l} \(H \nabla_i \a_{jl} - \a_{jl} \nabla_i H\)^2 \nn
\\ & \geq \frac1{|A|^2} \sum_{i,j,l} \dot \sigma^i  \(H \nabla_i \a_{jl} - \a_{jl} \nabla_i H\)^2, \lb{auxP1}
\end{align}
since by convexity $0<k_j<|A|$ for all $j$.

We now use the property that each $\dot \sigma^i$ is positive in the interior of the positive cone. More precisely, let us set, for any $\eps \in (0,1/n]$
$$
{\cal K}_\eps:=\{ \frak k =(k_1,\dots,k_n) \in \re^n ~:~ \min_{1 \leq i \leq n} k_i \geq \eps(k_1+\dots+k_n) >0 \, \},
$$ 
$$
M_1(\eps)=\min \{ \dot \sigma^i (\frak k) ~:~  1 \leq i \leq n, \ {\frak k} \in {\cal K}_\eps, | {\frak k}|=1 \}.
$$
Observe that $\dot \sigma^i >0$ on ${\cal K}_\eps$ for all $i$, by Lemma \ref{lpolS} (a). Therefore $M_1(\eps)>0$, being the minimum of a finite family of positive smooth functions on a compact set. In addition, since the cone ${\cal K}_\eps$ becomes smaller as $\eps$ increases, $M_1(\eps)$ is an increasing function of $\eps$. By homogeneity, we conclude
$$
\dot \sigma^i( \frak k) \geq M_1(\eps) |\frak k|^{m\beta-1}, \qquad {\frak k} \in {\cal K}_\eps.
$$

Substituting this in \eqref{auxP1} and using Lemma \ref{leHuS} we obtain that, on a hypersurface satisfying $\alpha \geq \eps Hg$, we have
\begin{align}
\left|H \nabla \a - \nabla H \, \a\right|^2_{\dot \sigma,b} & \geq  
\frac{n-1}{2}M_1(\eps)  \eps^2  |A|^{m \beta -3} H^2 |\nabla A|^2.
 \lb{auxP2}
\end{align}

We now want to estimate from above the term $\big|\ddot \sigma(\nabla A, \nabla A)\big|$. 
Observe that the quantity $\ddot \sigma(\nabla A, \nabla A)$ is homogeneous of degree $m\beta-2$ in the curvatures and quadratic in $\nabla A$. It is smooth as long as the curvatures are all positive, while it may be in general not defined when one or more curvatures vanish. With an argument similar to the previous one, we see that, for any $\eps \in (0,1/n]$, there exists a constant $M_2(\eps)$ such that, at any point where $\alpha \geq \eps H g$, 
\bec
\big|\ddot \sigma(\nabla A, \nabla A)\big| \leq M_2(\eps)|A|^{m \beta-2} |\nabla A|^2. \lb{auxP3}
\eec
The constant $M_2(\eps)$ is decreasing in $\eps$, since it gives a bound from above.

To conclude, we show that $\left| b - \frac{n}{H} g^{-1}\right|$ is small if the second fundamental form is  pinched enough. Clearly, we have
$$
 \left| b - \frac{n}{H} g^{-1}\right| \leq \max \left\{ \sqrt n \left( \frac{1}{k_1} - \frac{n}{H}  \right) , \sqrt n \left(\frac{n}{H} - \frac{1}{k_n}\right) \right\}.
$$
If $\alpha \geq \eps Hg$ for some $\eps \in (0,1/n]$, then $k_1 \geq \eps H$ and $k_n \leq (1-(n-1)\eps)H$. It follows
$$
\frac{1}{k_1} - \frac{n}{H} = \frac{H-nk_1}{k_1 H} \leq \frac{1-n\eps}{\eps H},
$$
$$
\frac{n}{H} - \frac{1}{k_n} = \frac{nk_n-H}{k_n H} \leq \frac{(n-1)(1-n\eps)}{k_n} \leq  \frac{n(n-1)(1-n\eps)}{H}.
$$
Since $\eps \leq 1/n$, we deduce that
\bec \lb{auxP4}
 \left| b - \frac{n}{H} g^{-1}\right| \leq (1-n\eps)\frac{n^{3/2}(n-1)}{H}.
\eec
Plugging \eqref{auxP2}, \eqref{auxP3} and \eqref{auxP4} into \eqref{auxP0}, we obtain
\begin{align}
D_+ \frak q &\geq \frac{n-1}{2} \frak q \, |A|^{m \beta -3} |\nabla A|^2 \(M_1(\eps)\eps^2 - 2 n^{3/2} (1-n\eps)\frac{|A|}{H} M_2(\eps) \) \nonumber
\\ & \geq   \frac{n-1}{2} \frak q \, |A|^{m \beta -3} |\nabla A|^2 \(M_1(\eps)\eps^2 - 2 n^{3/2} (1-n\eps) M_2(\eps) \). \label{auxP5}
\end{align}

To apply the maximum principle, we need that $M_1(\eps)\eps^2 - 2 n^{3/2} (1-n\eps) M_2(\eps) \geq 0$ on our hypersurface. Since $M_1(\eps)$ is increasing and $M_2(\eps)$ is decreasing, such a quantity is a strictly increasing function of $\eps$; in addition, it is negative for $\eps$ close to zero and positive for $\eps$ close to $1/n$. The optimal condition is obtained if we fix $\eps \in (0, 1/n)$ to be the unique value such that
\bec \label{eqeps}
M_1(\eps)\eps^2 - 2 n^{3/2} (1-n\eps) M_2(\eps) = 0.
\eec
By Lemma \ref{lem.const} there exists a constant $C_\fp \in \,(0,1/n^n)\,$ such that $K/H^n > C_\fp$ implies $\a > \eps Hg$, with $\eps$ given by \eqref{eqeps}. Then, if $K/H^n > C_\fp$ everywhere on our hypersurface, we have $D_+ \frak q \geq 0$ by \eqref{auxP5}. By the maximum principle, this proves that, for any $C > C_\fp$, the property $K/H^n \geq C$ is invariant under the flow.
\edem

\bnp \lb{const}
The quantities $M_1(\eps), M_2(\eps)$ introduced in the previous proof depend, in addition to $\eps$, only on the dimension $n$ and on the parameters $m,\beta$ which appear in the definition \eqref{def_sig} of the speed. Therefore the value of $\eps$ such that \eqref{eqeps} holds, and the constant $C_{\fp}$, only depend on $n,m,\beta$.
\enp

Theorem \ref{pinching} states that inequality $K/H^n \geq C_{\frak p}$ holds for all $t \in [0,T)\,$; in addition, by the definition of $C_{\frak p}$, we have that 
\bec
 k_i \geq \eps H \quad \text{on } \quad M \times [0, T) \quad \text{ for each $i$}  \lb{i+ii}
\eec
with $\eps$ given by (\ref{eqeps}). In particular, the solution is convex for all $t$ and therefore satisfies
  \bec
  k_j \leq H \quad \text{ on } \quad M \times [0, T) \quad \text{ for each $j$}. \lb{kmenH}
  \eec
Actually,  property (\ref{i+ii}) implies the sharper inequality $k_j \leq (1-(n-1)\eps)H$, but for our purposes it will suffice to use the simpler one (\ref{kmenH}).

Another consequence of the theorem is  a uniform double side bound for the inradius and outer radius of the evolving hypersurfaces.
\bco
\lb{db_in_kb} 
Under the assumptions of Theorem \ref{pinching}, there are constants $c_i = c_i(n, m, \beta, V)$, $i \in \{1, 2\}$,  such that 
\bec \lb{db_inr}
c_1 \leq \rho_t \leq D_t \leq c_2\qquad \text{ for every } \quad t\in [0,T)\,.
\eec
\eco

\bdem
If we denote by $\Omega_t$ the region enclosed by $M_t$, we have by the volume-preserving property of the flow and the definitions of $\rho_t,D_t$,
\begin{align*}
\omega_{n+1}  \rho_t^{n + 1}  \leq  \vle(\Omega_t) \leq \omega_{n+1}  D_t^{n + 1} ,
\end{align*}
where $\omega_{n+1}$ denotes the volume of the unit ball in $\re^{n+1}$. Since $\vle(\Omega_t) \equiv V$, we deduce
\bec \lb{radii}\rho_t \leq \(\frac{V}{\omega_{n+1}}\)^{\frac1{n + 1}} \leq D_t.\eec

On the other hand,  by (\ref{i+ii}) and (\ref{kmenH}), we have that $k_n \leq \eps^{-1} k_1$ on $M_t$.
Consequently, Lemma \ref{lem_An} implies that $D_t \leq B \rho_t$ for some constant $B=B(\eps,n)$. Combining this with \eqref{radii}, we reach the conclusion.
\edem

\section{Uniform bound for the velocity of the flow}

In this section we show that the pinching estimate implies a uniform bound from above for the speed of the flow and for the curvature of the hypersurface. Throughout the section, we assume that the flow satisfies the assumptions of Theorem \ref{pinching}. As usual, we denote by $\Omega_t \subset \re^{n+1}$ the region enclosed by $M_t$. Following the procedure of \cite{MC04,Mc3}, we first prove a result about the existence of a ball with fixed center enclosed by our hypersurface on a suitable time interval.

\bp \lb{comp_pple} 
Given any $\bar t \in [0,T)$, let $\bar x \in \Omega_{\bar t}$ be such that $B(\bar x,\bar \rho) \subset \Omega_{\bar t}$, where $\bar \rho:=\rho_{\bar t}$ is the inradius of $M_{\bar t}$. Then we have
\bec \lb{comp_p}
B(\bar x, \bar \rho/2) \subset \Omega_t \qquad \text{ for every } \quad t\in [\bar t, \min\{T,  \bar t+\tau\}) \,.
\eec
for some constant $\tau$ depending only on $n, \beta, m, V$.
\ep

\bdem
Given $\bar x, \bar \rho$ such that $B(\bar x,\bar \rho) \subset \Omega_{\bar t}$, let us denote by
$\frak X(\, \cdot \,, t) = X(\, \cdot \,, t) - \bar x$ the position vector field (with origin $\bar x$) on the evolving manifold $M_t$ and define $r(\, \cdot \,, t) := |\frak X(\, \cdot \,, t)|$. Notice that $r^2$ evolves under  \eqref{vpgcf} according to
$$\parcial{r^2}{t} = \parcial{}{t} \<\frak X, \frak X\> = 2 \<\frak X, \parcial{\frak X}{t}\> = 2 \, (h - \sigma) \<\frak X, N\>.$$
It is easy to check that 
$\Delta r^2 = - 2 H  \<\frak X, N\>   + 2 n$. 
 Thus we can write
\bec \lb{eq_r2}
\(\parcial{}{t} - \frac{\sigma}{H} \Delta\) r^2  = 2 h  \<\frak X, N\> - 2 n  \frac{\sigma}{H}
 \geq 2  \, h \<\frak X, N\> - 2 n \frac{H^{m \beta -1}}{n^{m \beta}}, 
\eec
where we have used part (d) of Lemma \ref{lpolS}. 

We set $t_1 = \inf\{t > \bar t \, / \, \bar x \not \in \Omega_t\}$ provided this set is nonempty; we set $t_1=T$ otherwise. We then have 
\bec \lb{bor_cond2}
t_1<T \quad \Longrightarrow \quad \bar x \in M_{t_1} = \partial \Omega_{t_1}. 
\eec
In addition, by convexity, we have $\<N, \frak X\> \geq 0$  on $M_t$ for all $t \in [\bar t, t_1)$. Then the evolution equation for $r$ on $[\bar t, t_1)$ becomes
\begin{align}
\(\parcial{}{t} - \frac{\sigma}{H} \Delta\) r & = \frac1{2 r} \(\parcial{}{t} - \frac{\sigma}{H} \Delta\) r^2  + \frac{1}{r} \frac{\sigma}{H} |\nabla r|^2\nn
\\ & \usb{eq_r2}{\geq} \frac1{r}\<\frak X, N\>  h - 	\frac1{r} \(\frac{H}{n}\)^{m \beta -1} \geq - \frac1{r} \(\frac{H}{n}\)^{m \beta -1}. \lb{cpw_fin} 
 \end{align}

Next, we define $\frak r(t):= \min_M r(\, \cdot \,, t)$ for any $t \in [\bar t, t_1)$ and set
$Y(t):= \{p \in M \, / \, r(p, t) = \frak r(t)\}.$
Now \eqref{cpw_fin} allows us to deduce
$$D_+ \frak r(t) \geq \inf_{Y(t)} \parcial{r}{t} \geq  \inf_{Y(t)} \left[ - \frac1{r} \(\frac{H}{n}\)^{m \beta -1}\right].$$
Observe that at any point where the minimum is attained the hypersurface is tangent to an inball of radius $\frak r(t)$, which implies that $H \leq n/\frak r(t)$ on any point of $Y(t).$
 As, in addition, $r(\, \cdot \,, t)$ is constantly equal to $\frak r(t)$ on $Y(t)$, we conclude 
$$D_+ \frak r(t) \geq - \frac1{\frak r(t)^{m \beta}}.$$
Using a standard comparison principle (cf. \cite [\S 10.3]{ChII}), we obtain
\bec \lb{IVP_R}
\frak r(t) \geq R(t), \qquad  \text{for } \quad t \in [\bar t, t_1),
\eec
being $R$ the solution to $R' = - R^{- m \beta}$ with the initial value $R(\bar t) = \bar \rho$. 
Notice that $R(t)$ can be regarded as the radius of a geodesic sphere $\partial B(\bar x, R(t))$ contracting under the ordinary $H_m^\beta$-flow  (that is, \eqref{vpgcf} setting $h = 0$) with initial condition $\partial B(\bar x, \bar \rho)$.

Setting $d:= m \beta + 1$, from the explicit expression of $R$, we conclude that
\begin{align*} 
R(t)  =\(\bar \rho^{d} - d \,(t - \bar t)\)^{\frac1{d}}\geq  \frac{\bar\rho}{2} \qquad \text{if and only if} \qquad  t - \bar t \leq 
 \frac{2^{d} - 1}{d} \(\frac{\bar \rho}{2}\)^{d}.
\end{align*}
and thus
\bec \lb{cotaR}
R(t) \geq \frac{\bar\rho}{2} \qquad \text{ if } \qquad t - \bar t \leq \frac{2^{d} - 1}{d}\(\frac{c_1}{2}\)^{d} =: \tau,
\eec
where $c_1$ is the constant coming from the lower bound of Corollary \ref{db_in_kb}. 

To complete the proof, assume that $t_1 <  \min\{\tau + \bar t, T\}$. By \eqref{IVP_R},  we have
$$r(q, t_1 - \eta) = |X(q, t_1 - \eta) - \bar x| \geq R(t_1 - \eta) \quad \text{ for all } q \in M_{t_1 - \eta} ,\ \eta \in \,(0,t_1-\bar t].$$
Hence 
$$0 \us{\eqref{bor_cond2}}{=} \dist(M_{t_1}, \bar x) = \lim_{\eta \to 0^+} \dist(M_{t_1 - \eta}, \bar x) \geq R(t_1) \us{\eqref{cotaR}}{\geq} \bar \rho/2 > 0,$$
which is a contradiction. In conclusion, $t_1 \geq \min\{\tau + \bar t, T\}$ which, together with   \eqref{IVP_R} and \eqref{cotaR}, leads to
$$\frak r(t) \geq \bar \rho/2 \qquad \text{on } \quad [\bar t, \min\{\bar t + \tau, T\}),$$
which proves our assertion.
\edem

The above result allows us to obtain a uniform bound on the speed of the flow using a technique that was first introduced by Tso \cite{Tso}.

\bp \lb{analogSch}
  There exists a constant $C =  C(n, m, \beta , M_0)$ such that
\bec \lb{bo_Hm}
H_m (\cdot, t) \leq C \qquad \text{ for all } \quad t\in [0,T).
\eec
Moreover, the upper bounds 
\bec \lb{cot_h_kb}
h(t) \leq C^\beta \qquad \text{ and }\qquad H(\cdot, t) \leq C_1 =C_1(C, \eps)
\eec
hold for every $t \in [0,T)$.
\ep

\bdem
For any fixed $\bar t \in [0,T)$, let $\bar x$ and $\bar \rho$ be as in Proposition \ref{comp_pple}. For brevity, let us set $u=\<\frak X, N\>$, where $\frak X$ is the position vector field with origin $\bar x$.
Using \eqref{comp_p} and the convexity of the $M_t$'s, and choosing $c:= \bar \rho/4$, we achieve
\bec 
u - c \geq \frac{\bar \rho}{2} - \frac{\bar \rho}{4} = c > 0  \qquad \text{on } [\bar t, \min\{\bar t + \tau, T\}),  \lb{u_c_pos}
\eec
which ensures  that the function $W \ds = \frac{\sigma}{u - c}$ is well-defined. A routine computation (as in \cite{MC04}, \cite{CaMi1}) using \eqref{evF} and \eqref{fal_sopev} gives
\begin{align*}
\(\parcial{}{t} - \Delta_{\dot \sigma}\) W = 2\frac{\<\nabla W, \nabla u\>_{\dot \sigma}}{u - c} + (m \beta + 1) W^2 - h \frac{W + \tr_{\dot \sigma}(\a A)}{u - c} - c \, W \frac{\tr_{\dot \sigma}(\a A)}{u - c}.   
\end{align*}

As strict convexity holds for each $M_t$, we have that $h$, $\sigma$ and $k_i$ are all positive, which together with \eqref{u_c_pos} allow us to disregard the term containing $h$. In this way, we obtain the following estimate:
\bec \lb{eqW}
\(\parcial{}{t} - \Delta_{\dot \sigma}\) W \leq \frac2{u - c}\<\nabla W, \nabla u\>_{\dot \sigma} + (m \beta + 1) W^2  - c \, m \beta \, \eps \, H W^2
\eec
on $[\bar t, \min\{\bar t + \tau, T\})$, where we have also used that
$$
 \tr_{\dot \sigma}(\a A) =  \dot \sigma^i k_i^2 \us{\eqref{i+ii}}{\geq} \eps H \dot \sigma^ i k_i = \eps \, m \beta \, \sigma  
$$
since $\sigma$ is homogeneous of degree $m\beta$. 

Using \eqref{u_c_pos} and Lemma \ref{lpolS} (d), we have
$$
W \leq \frac{\sigma}{c} \leq \frac 1c \left( \frac Hn \right)^{m\beta}, 
$$
which, together with \eqref{eqW}, yields
\bec \lb{eqW2}
\(\parcial{}{t} - \Delta_{\dot \sigma}\) W \leq \frac2{u - c}\<\nabla W, \nabla u\>_{\dot \sigma} + 
\left(m \beta + 1  - c^{1+\frac{1}{m\beta}} \,  m n \beta  \eps \, W^\frac 1{m\beta} \right) W^2.
\eec
From \eqref{eqW2},  using a maximum principle argument as in Corollary 4.5 in \cite{MC04}, we obtain \eqref{bo_Hm}.

The bound  \eqref{bo_Hm} on $H_m$ implies that $h \leq C^{\beta}$ by the definition of $h$.  Next, by homogeneity and the inequality (c) in Lemma \ref{lpolS}, we have 
$$m  \beta\, \sigma =  \dot \sigma^i k_i \us{\eqref{i+ii}}{\geq} \eps H \tr(\dot \sigma) \geq \eps  H \, m \beta\sigma^{1 - 1/m\beta}.$$
Then
\begin{align}
H & \leq \frac1{\eps} \sigma ^{1/m\beta} = \frac1{\eps} H_m^{1/m}, \lb{HloF}
\end{align}
and so $H \leq \eps^{-1} C^{1/m} =:C_1$ again by application of  \eqref{bo_Hm}.
\edem

Since our hypersurfaces are convex, the bound on $H$ we have just obtained implies a bound on all principal curvatures. As a consequence, one can prove that $M_t$ can be locally written as a graph with uniformly bounded $C^2$ norm (see e.g. the proof of Lemma 3.4 in \cite{Sch2}):

\bco \lb{graph_u} There exist $r,\eta>0$ (depending only on $\max H$) with the following property.
Given any $(\bar p, \bar t) \in M \times \,(0, T)$, there is a neighborhood $\mc U$ of the point $\bar x:= X(\bar p, \bar t)$ such that $M_t \cap \mc U$ coincides with the graph of a smooth function 
$$u: B_r \times J \flecha \re, \qquad \mbox{ for all } t \in J.$$
Here $B_r \subset T_{\bar p} M_{\bar t}$ is the ball of radius $r$ centered at $\bar x$ in the hyperplane tangent to $ M_{\bar t}$ at $\bar x$, and $J$ is the time interval $J=\(\max\{\bar t-\eta,0\}, \min\{ \bar t + \eta, T\}\)$. In addition, the $C^2$ norm of $u$ is uniformly bounded (by a constant depending only on $\max H$).
\eco
 
For future reference, it  is useful to recall here some basic formulae  relating quantities over $M_t$ with the function $u$ giving the local parametrization of $M_t$ as a graph over a hyperplane. Using $D_i$ to denote the derivatives with respect to these local coordinates, and choosing as positive normal the one which points below, we have (cf. \cite{Ur2}, \cite{Eck})
\bec \lb{met_u}
g_{ij} = \delta_{ij} + D_i u \, D_j u, \quad  \quad \quad g^{ij} = \delta^{ij} - \frac{D^i u \, D^j u}{1 + |D u|^2},
\eec
and
\bec \lb{a_u}
\a_{ij} = \frac{D_{ij} u}{(1 + |Du|^2)^{1/2}}.
\eec
In addition, the Christoffel symbols have the expression:
\bec \lb{invg_u}
\Gamma_{ij}^k = \(\delta^{kl} - \frac{D^k u \, D^l u}{1 + |D u|^2}\) D_{ij} u \, D_l u.
\eec

\section{Long time existence} \lb{lte_gcf}

In this section we prove that the solution of problem \eqref{vpgcf}--\eqref{inval} exists for all times if the initial immersion satisfies the pinching condition \eqref{in_pinch}. Let us denote by $[0,T_{\max})$ the maximal interval of existence of the solution, with $T_{\max} \in (0,+\infty]$. As usual, we shall prove that $T_{\max}=+\infty$ by a contradiction argument, showing that, if $T_{\max}$ is finite, then the solution can be continued for some times larger than $T_{\max}$. To this purpose, we need suitable estimates on the solution on any finite interval which guarantee that $M_t$ converges to a smooth limiting hypersurface as $t \to T_{\max}$. In addition, we have to show that the solution remains uniformly convex on $[0,T_{\max})$, to ensure that the limit still satisfies the parabolicity assumption.

We first consider the issue of the preservation of convexity. Observe that Theorem \ref{pinching} implies the uniform convexity of our hypersurfaces; however, that result was obtained under the additional apriori assumption that $H>0$. Such an assumption is certainly valid for $t$ close to zero, by \eqref{in_pinch}, but there may be some positive time at which both $\min K$ and $\min H$ go to zero, at a rate such that $K/H^n$ remains bounded. To exclude such a behavior, we have to complement Theorem \ref{pinching} with a lower bound on $H$ on any finite time interval, which is given by the next lemma.

\bl \lb{lobo_K}
Under the hypotheses of Theorem \ref{pinching}, there exist $C_2,C_3>0$ depending on $n, m, \beta, M_0$ such that 
\bec \lb{Hm_lo}
\min_{M_t} H \geq C_2 e^{ -C_3 t } \qquad \forall  t \in [0, T).
\eec
\el

\bdem
Since we are under the hypotheses of Theorem \ref{pinching}, we know that our hypersurfaces $M_t$ are convex for every $t \in [0,T)$. Here we want to show that $H$ satisfies the lower bound \eqref{Hm_lo}, which in particular implies that $\min H$ cannot tend to zero as $t \to T$.

We first derive a bound from below for $\sigma$. At a point where the minimum $\frak s(t) := \min_M \sigma(\, \cdot \,, t)$ is attained,  we deduce from \eqref{evF}
$$D_+ \frak s \geq  (\sigma - h)\tr_{\dot \sigma} (\a A)  \geq  - h \, H \, m \beta \, \frak s \geq - C^\beta C_1 m \beta \, \frak s =: - \gorro C \frak s,$$
where we have used that $\sigma \, \tr_{\dot \sigma} (\a A) \geq 0$ (by convexity) 
and the computation
\bec
\tr_{\dot \sigma} (\a A) = \dot \sigma^i k_i^2 \us{\eqref{kmenH}}{\leq} H \dot \sigma^i k_i = H m \beta \sigma, \lb{traA_up}
\eec
which follows  by homogeneity of $\sigma$. Then, a scalar maximum principle leads us to
\bec \lb{lo_K_T}
\frak s(t) \geq \frak s(0) \, e^{-\gorro C\, t}.
\eec
By Lemma \ref{lpolS} (d) we conclude that
$$
\min_{M_t} H \geq n \, {\frak s(t)}^{1 / m\beta} \geq n  (\min_{M_0} \sigma )^{1 / m\beta} e^{-(\gorro C /m\beta)t},
$$
which implies the assertion.\edem

\bco \lb{pos_H}
Let $X:M \times [0,T_{\max}) \, \to \re^{n+1}$ be the solution of \eqref{vpgcf} with an initial value which satisfies the pinching condition \eqref{in_pinch}. Then, the hypersurfaces $M_t$ are uniformly convex on any finite time interval; that is, for any $T < +\infty$, $T \leq T_{\max}$, we have
$$
\inf_{M \times [0,T)} k_i >0, \qquad \forall i=1,\dots,n.
$$
Therefore, Theorem \ref{pinching} is valid also without the hypothesis that $H>0$ for $t \in (0,T)$. The same holds for the other results that have been obtained until here under the same assumptions of Theorem \ref{pinching}.
\eco
\bdem
The pinching condition \eqref{in_pinch} implies that the initial immersion is uniformly convex. Therefore for $T$ small enough the corollary is true by continuity. If the corollary is not true in general, there exists a smallest time $T_0>0$ such that $\inf_{M \times [0,T_0)} k_i =0$ for some curvature $k_i$. But then $H>0$ for $t \in [0,T_0)$; therefore we can apply Theorem \ref{pinching} and Lemma \ref{lobo_K} on this time interval to conclude that
$$
\inf_{M \times [0,T_0)} k_i  \geq \eps \inf_{M \times [0,T_0)} H \geq \eps \, C_2 \, e^{-C_3T_0}>0,
$$
which is a contradiction.
\edem

We now have to show that, on any finite time interval $[0,T)$, the hypersurfaces $M_t$ satisfy uniform estimates on the curvature and all its derivatives, in order to guarantee that, if $T_{\max}$ is finite, then they converge to a smooth limit as $t \to T_{\max}$. This part is technically more complicate. It will be convenient to represent locally $M_t$ as the graph of a function $u$ and obtain uniform estimates on the derivatives of $u$. We already have a $C^2$ estimate on $u$  coming form the curvature bounds obtained in the previous section. It is well known that the crucial step is to derive uniform $C^{2,\alpha}$ estimates for some $\alpha>0$, since such a property allows us to deduce all the higher order estimates by using standard linearization and bootstrap techniques. However, the parabolic equation solved by our hypersurface is fully nonlinear and in this case $C^{2,\alpha}$ estimates are known in general only for operators with concave dependence on the second derivatives \cite{K,Lieb}. In our case, instead, the operator is not concave, as it can be clearly seen from the property that it is homogeneous with degree larger than one.

However, our operator is not too far from being concave because it is, roughly speaking, a power of a concave operator. In fact, the speed $\sigma$ can be written as the $m\beta$-th power of $H_m^{1/m}$ which is a concave function, as we have recalled in Lemma \ref{lpolS} (b). This allows us to use a technique inspired by the work of B. Andrews \cite{A5}, D. Tsai \cite{Tsai} and J. McCoy \cite{Mc3}, where the crucial step is to prove space regularity at a fixed time by applying Theorem \ref{Caf-pert}. We prove this step separately in the following lemma.

\bl \lb{appl_caff}
Let $M \subset \re^{n+1}$ be an embedded hypersurface satisfying at every point $H_0<H<H_1$, $k_1 \geq \eps H$ for given positive constants $H_0,H_1,\eps$. Given any $p \in M$, let $u$ be a local graph representation of $M$ over a ball $B_r \subset T_pM$ (as in Corollary \ref{graph_u} applied at a fixed time). Then $u$ satisfies
$$||u||_{C^{2,\alpha}(B_r)} \leq C(1+||\sigma||_{C^{\alpha}(B_r)})$$
for some $C>0$ and $0<\alpha<1$ depending only on $n,H_0,H_1,\eps$ and the parameters $\beta,m$ in the definition of $\sigma$. 
\el

\bdem
The hypotheses imply that the principal curvatures of $M$ at every point are contained between two fixed positive constants. We know from Corollary \ref{graph_u} that $M$ can be written locally as the graph of a function $u$ with $||u||_{C^2}$ bounded in terms of $H_1$.  To show that $u$ satisfies a $C^{2,\alpha}$ estimate as well, we want to apply Theorem \ref{Caf-pert} to the function $u$, choosing as $F$ the function expressing $\sigma$ in the graph representation of $M$. 

Let us observe that Theorem \ref{Caf-pert} is stated for operators defined for arbitrary values of the hessian, while $\sigma$ in general is well-defined and elliptic only if $D^2u$ is positive definite. However, this is not a substantial difficulty since we are assuming a priori uniform convexity of our hypersurface. 
An easy way to handle this problem is replacing $\sigma$ by another function $\tilde \sigma$, which is defined everywhere and coincides with $\sigma$ on a set containing the possible values of the curvatures of $M$. To this purpose, let us first define
$$
{\cal C}=\{ \frak k=(k_1,\dots,k_n) ~:~ H_0 \leq H(\frak k) \leq H_1, \quad \min_{1 \leq i \leq n}k_i \geq \eps H(\frak k) \},
$$
which is a compact symmetric subset of the positive cone $\Gamma_+$. For simplicity of notation, let us set $\phi:=H_m^{1/m}$. We know from Lemma \ref{lpolS} that $\phi$ is homogeneous of degree one and is concave in $\Gamma_+$. This implies that, given any $\frak h, \frak k \in \Gamma_+$, we have
\bec \lb{supd}
\phi(\frak k) \leq \phi({\frak h}) + D\phi ({\frak h}) \cdot ({\frak k}-{\frak h}) =
D\phi ({\frak h}) \cdot {\frak k},
\eec
since $D\phi ({\frak h}) \cdot {\frak h} = \phi({\frak h})$ by Euler's theorem on homogeneous functions. Let us now define
\bec \lb{defBE}
\tilde \phi ({\frak k})=\min_{{\frak h} \in \cal C}  D\phi ({\frak h}) \cdot {\frak k}, \qquad \frak k \in \re^n.
\eec
The minimum is attained for any ${\frak k} \in \re^n$ because $\cal C$ is a compact set. The function $\tilde \phi$ is (positively) homogeneous of degree one, by definition, and is concave, since it is the minimum of linear functions.  Using (\ref{supd}) and the property that $D\phi ({\frak h}) \cdot {\frak h} = \phi({\frak h})$, we see that $\tilde \phi$ coincides with $\phi$ on $\cal C$; in addition, $\tilde \phi$ is symmetric, because $\phi$ and $\cal C$ are both symmetric. Let us also observe that, since $\partial \phi / \partial k_i (\frak k) >0$ for any $\frak k \in {\cal C}$, and $\cal C$ is compact, there exist $m_2>m_1>0$ such that
$$
m_1 \leq \partial \phi / \partial k_i (\frak k) \leq m_2, \qquad i=1,\dots,n, \ \frak k \in {\cal C}.
$$
This implies easily, using the definition of $\tilde \phi$, that
\begin{equation}\label{ell.phi}
m_1 |\frak l| \leq  \tilde \phi ({\frak k}+{\frak l}) - \tilde \phi ({\frak k}) \leq \sqrt n \,  m_2  |\frak l|, \mbox{ for all }\frak k, \frak l \in \re^n, \ {\frak l} \geq 0,
\end{equation}
where ${\frak l} \geq 0$ means that all components of ${\frak l}$ are nonnegative.

Now let $u:B_r \to \re$ be a local graph representation of $M$ over its tangent plane at a given point, and let the coordinate in $B_r$ be denoted by $x=(x_1,\dots,x_n)$. Let us consider the function $\tilde \phi({\frak k}(x))$, where ${\frak k}(x)$ are the principal curvatures of $M$ at the point $(x,u(x))$. Since ${\frak k}(x)$ are the eigenvalues of a matrix depending on $Du,D^2u$ (see \eqref{met_u}, \eqref{a_u}), we have that $\tilde \phi({\frak k}(x))$ can be expressed as $\tilde \Phi(Du(x),D^2u(x))$ for a suitable function $\tilde \Phi = \tilde \Phi(p,A)$, with $(p,A) \in \re^n \times {\cal S}$, ${\cal S}$ being the set of symmetric $n \times n$ matrices.
The dependence of $\tilde \Phi$ on $A$ is related to the dependence of $\tilde \phi$ on $\frak k$.
In fact, it is well known (see e.g. \cite{CNS4, Ur2,Ger}) that the concavity of $\tilde \phi$ with respect to $\frak k$ implies the concavity of $\tilde \Phi$ with respect to $D^2 u$, and that (\ref{ell.phi}) implies the ellipticity condition (\ref{FNL-ell}) for $\tilde \Phi$. In addition, $\tilde \Phi$ is homogeneous of degree one with respect to $D^2u$.

Now let us further set $F(x,A):=\tilde \Phi(Du(x),A)$ and $f(x)=\tilde \Phi(Du(x),D^2u(x))$. Then, $u$ can be regarded as a solution of the equation
$$
F(x,D^2u(x))=f(x), \qquad x \in B_r.
$$
The above remarks show that all hypotheses of Theorem \ref{Caf-pert} are satisfied and that there exists $\a \in (0,1)$ such that 
$$\|u\|_{C^{2+\a}(B_{r/2} )} \leq C(1+||f||_{C^\alpha(B_r)}),$$
where $C$ depends on $H_0,H_1,\eps,\beta,n$.

By our assumptions, ${\frak k}(x)$ belongs for every $x$ to the set $\cal C$ where $\tilde \phi$ and $\phi$ coincide. Therefore, $f$ coincides with $H_m^{1/m}=\sigma^{1/\beta m}$ evaluated at ${\frak k}(x)$. Observe that our assumptions on the curvatures imply that $\sigma$ is contained between two positive values depending only on $H_0,H_1,\eps,n,\beta$. Therefore $||\sigma^{1/\beta m}||_{C^\alpha}$ is estimated by $||\sigma||_{C^\alpha}$ times a constant depending only on these quantities. This completes the proof. 
\edem

The function $\tilde \phi$ was used also in \cite{Mc3}, where it was called {\em Bellman's extension} of $\phi$. However, we have slightly changed the definition because of the following problem. In \cite{Mc3} the function $\tilde \phi$ is defined as in \eqref{defBE} with $\cal C$ replaced by the positive cone. However, functions of the form $H_m^{1/m}$ are in general not differentiable on the boundary of the positive cone. In fact, $H_m$ vanishes at some points of the boundary and the gradient of $H_m^{1/m}$ becomes unbounded as such points are approached. Consequently, if $\cal C$ is replaced by the positive cone the minimum in \eqref{defBE} is in general not attained and $\tilde \phi$ is equal to $-\infty$ in some regions.

\bt \lb{teor_reg}
Let $M_t$ be a solution of \eqref{vpgcf}, defined for $t \in [0,T_{\max})$, with initial condition satisfying \eqref{in_pinch}. Then, for any $0<\delta<T \leq T_{max}$, with $T<+\infty$ and any integer $l \geq 0$, there exists $\mc C=\mc C(n, m, \beta, l,\delta,T, M_0)$ such that
$$
\sup_{M \times [\delta,T) } |\nabla^l A|  \leq \mc C.
$$
\et

\bdem 
The case $l=0$ follows from the curvature bounds of the previous sections. To consider a general $l$,  we write locally $M_t$ as the graph of a function $u$ as in Corollary \ref{graph_u}. Then $u$ satisfies the equation
\bec \lb{graph}
\frac{\partial u}{\partial t} = \sqrt{1+|Du|^2} (\sigma-h(t)).
\eec
The function $\sigma$ depends on the Weingarten operator associated to $M_t$; hence, in the coordinate system under consideration, it is a function of $D^2 u$ and $Du$. The right hand side is a fully nonlinear operator; as we have already observed, it is an elliptic operator, thanks to property (a) of Lemma \ref{lpolS}. The higher order regularity does not follow by the general theory of Krylov and Safonov \cite{K,Lieb} because the operator is not a concave function of $D^2 u$. We use instead the procedure of \cite{A5,Tsai,Mc3}, which consists of proving first regularity in space at a fixed time and then regularity in time.

We start by deriving a $C^{ \a}$-estimate for $\sigma$.  To do so, let us denote by $D_i$ the derivatives with respect to the local coordinates in the graph representation of $M_t$. Then equation \eqref{ev_sig2} can be written as
\begin{align} 
\parcial{\sigma}{t} &= a^{ij} D_i D_j \sigma + b^i D_i \sigma + e \, \sigma, 
\qquad (x,t) \in B_r \times J, 
\lb{Kb_KS}
\end{align} 
with $B_r,J$ as in Corollary \ref{graph_u}, and the coefficients given by
$$a^{ij} = \beta \, H_m^{\beta - 1}\, c^{ij}, \quad b^i = \beta\, H_m^{\beta - 1} \, c^{lj} \Gamma_{lj}^i \quad \text{and } \quad e = \beta\,(\sigma - h) H_m^{-1}\tr_c (\a A).$$
Since we have uniform bounds on the curvatures both from above and from below on any finite time interval, equation \eqref{Kb_KS} is uniformly parabolic with uniformly bounded coefficients. Then we are in position to apply Theorem \ref{KS_t} and obtain
\bec
\|\sigma\|_{C^\a(B_{r'} \times J^*)} \leq \mc C \|\sigma\|_{C^0(M \times [0,T))} \leq  \mc C' \lb{bo_sig}
\eec
for any $B_{r'} \subset B_r$, and a suitable $\a \in \, (0, 1)\,$, where $J^*$ is the interval $J$ minus its initial part of length $\delta$. 
By covering $M_t$ with graphs over balls of radius $r'$, we obtain an estimate on $||\sigma||_{C^\a(M \times [\delta, T))}$. Thus, if we fix any time $t \in [\delta,T)$, we can apply Lemma \ref{appl_caff} to $M_t$ and conclude that any graph representation $u(\cdot,t)$ of $M$ satisfies a uniform $C^{2,\alpha}$ estimate in space.

From this estimate on $u(\cdot,t)$ for any fixed $t$, we can deduce a $C^{2,\alpha}$ estimate for $u$ with respect to both space and time following the procedure of  \cite[\S 3.3, 3.4]{A5} or \cite[Th. 2.4]{Tsai} to equation \eqref{graph}. Once $C^{2,\alpha}$ regularity is established, standard parabolic theory yields uniform $C^k$ estimates for any integer $k > 2$. Since this holds for any graph representation of $M_t$, we deduce that any derivative $|\nabla^l A|$ of the Weingarten map is uniformly bounded on any finite time interval. 
\edem

With the above result, standard continuation techniques yield (see e.g. Theorem 8.1 in \cite{Hu84}):

\bt \lb{lte_kb} Let $M$ be a closed $n$-dimensional smooth manifold and $X: M \fle \re^{n + 1}$ be an immersion pinched in the sense of \eqref{in_pinch}. If $M_t = X_t(M)$ is the solution to \eqref{vpgcf} with initial condition $X_0 = X$, then  $M_t$ exists on $[0, \infty)$.
\et

The above result, together with Theorem \ref{pinching} and Corollary \ref{pos_H}, completes the proof of part (a) and (b) of our main Theorem \ref{t_main}. It remains to prove the convergence to a sphere, which will be done in the next section.

\section{Convergence to a round sphere}

Henceforth, $M_t = \partial \Omega_t$ stands for a solution of the flow  \eqref{vpgcf} on $[0, \infty)$, and with initial condition satisfying \eqref{in_pinch}. By our previous results, we know that the solution is smooth and uniformly convex on any finite time interval, but we cannot completely control its behavior as $t \to \infty$. In fact, since the lower bound for $H$ given by Lemma \ref{lobo_K} depends on $t$, we cannot exclude at this stage that $\min_{M_t} H \to 0$ as $t \to \infty$ so that uniform convexity is lost. Since the regularity estimates depend on uniform convexity, some additional argument is needed to ensure the existence of a smooth limit as $t \to \infty$.

Our strategy will consist of showing first that, if a smooth limit exists, it has to be a round sphere; the existence of the limit will be proved afterwards. To address the first step, we consider again the quotient $K/H^n$ that we have used for the proof of the preservation of pinching, and prove that it converges uniformly to $1/n^n$, which is the value assumed on a sphere.

\subsection{Asymptotic behavior of the quotient $\bm{K/H^n}$}

In order to analyze the behavior of $K/H^n$, we need to ensure that  the global term $h(t)$ is bounded from below by a positive constant. 

\bl \lb{loh}  There exists a constant $h_0 = h_0(n, m, \beta, V) > 0$ such that $h(t) \geq h_0$ for all $t \geq 0$. 
\el

\bdem
As our evolving manifolds are convex, we can use Lemma \ref{Alex_Fen} to estimate
\begin{align*}
\(C_n \, V^{n-1}\)^{\frac1{n + 1}} & \leq \int_M H \, d \mu_t \us{\eqref{HloF}}{\leq} \frac1{\eps} \int_M \(H_m^\beta\)^{\frac1{m \beta}} \, d \mu_t  \\
& \leq \frac1{\eps} \(\int_M H_m^\beta \, d \mu_t\)^{\frac1{m \beta}} |M_t|^{1-{\frac1{m \beta}}},
\end{align*}
where $V = \vle(\Omega_0)$ and we have applied a H\"older inequality.
Hence
\bec \lb{esti_h_Mt}
 h = \frac1{|M_t|} \int_M H_m^\beta \, d \mu_t \geq \eps^{m\beta} \(C_n \, V^{n-1}\)^{\frac{m \beta}{n +1}} \,   |M_t|^{- m \beta}.
 \eec

Now, by \eqref{db_inr}, we know that the area of $M_t$ is not greater than the area of a sphere of radius $c_2$. The use of this in \eqref{esti_h_Mt} gives the desired lower bound for $h$.
\edem

As we mentioned above, we are not yet able at this stage to exclude that $ \min_{M_t} H \to 0$ asymptotically. However, we show in the next lemma that $\min H$ does not decay too fast, and this will be enough for our purposes.

\bl \lb{0.19} We have
\bec \lb{int_H_inf}
\int_{0}^{\infty} H_{\min}(t) dt =+\infty,
\eec
where we have set $H_{\min}(t) = \min_M H(\, \cdot \,, t)$. 
\el

\bdem
Let us first estimate  $\frak s(t) := \min_M \sigma(\, \cdot \,, t)$. At a point where the minimum is attained, from \eqref{evF} we get 
\begin{align*}
{D_+ \frak s} & \geq  \(\sigma - h\) \tr_{\dot \sigma} (\a A).
\end{align*} 
As $\sigma \, \tr_{\dot \sigma} (\a A)$ is always nonnegative because $M_t$ is convex, we obtain
$$D_+ \frak s \geq  - h \, \tr_{\dot \sigma} (\a A) \us{\eqref{traA_up}}{\geq} -C^\beta \, m \beta \, H \, \frak s \us{\eqref{HloF}}{\geq} - \frac{m \beta}{\eps} C^\beta \frak s^{1 + \frac1{m\beta}} =: - \gorro C \frak s^{1 + \frac1{m\beta}}, $$ being $C^\beta$ the upper bound for $h$ coming from Proposition \ref{analogSch}. Now the application of a maximum principle to the above inequality yields
$$\frak s \geq  \(C_0 + \gorro C \, \frac{t}{m \beta}\)^{-m \beta}, \qquad \text{with} \quad C_0 := \frak s(0)^{-\frac1{m\beta}},$$
which, together with part (d) of Lemma \ref{lpolS}, leads to
$$\int_0^\infty H_{\min} (t) \, dt \geq \lim_{s \to \infty} \int_0^s n \,\frak s^{\frac1{m\beta}} (t) \, dt \geq \lim_{s \to \infty}  \int_0^s n \(C_0 + \gorro C \, \frac{t}{m \beta}\)^{-1} = \infty.$$
\edem

\bp  \lb{asym_pc}
The quotient $K/H^n$ converges to $1/n^n$ uniformly on $M$ as $t \to \infty$.
\ep

\bdem
We consider the evolution equation for  $q(\ccdot, t):= K/H^n$ given in \eqref{pin_vpgcf}.
First, applying Lemma \ref{loh},  Lemma  \ref{2.5} and inequality \eqref{in_pinch}, we estimate the term containing $h$ as follows:
\begin{align}
h \frac{q}{H}(|A|^2 n-H^2) & \geq  h_0 \, q\, H \,
\frac{|A|^2 n-H^2}{H^2} \geq  B \, H_{min}(t) \,  \(n^{-n} - q\), \lb{esti_th}
\end{align}
with $B := h_0 \,C_{\frak p} \, \delta$. Let us now define $\frak d(t) = \sup_M \(n^{-n} - q\)(\, \cdot \,, t)$. Then $\frak d$ is a locally Lipschitz continuous function and  satisfies 
$$D^+ \frak d(t) \leq \sup_{\mc M(t)} \parcial{}{t} \(n^{-n} - q\) = \sup_{\mc M(t)} \parcial{}{t}\(- q\),$$
where $\mc M(t)= \left\{p \in M \, / \, \frak d(t) = n^{-n} - q(p, t)\right\},$
 and $D^+$ stands for  the upper right Dini derivative (cf. \cite[\S 10.3]{ChII}). Using \eqref{pin_vpgcf}, we obtain
\begin{align*}
D^+ \frak d(t) & \us{\eqref{esti_th}}{\leq} \sup_{\mc M(t)} \left[- B \, H_{\min}(t) \(n^{-n} - q\)\right] = - B\, H_{\min}(t) \, \frak d(t),
\end{align*}
since the sum of the gradient terms in the second row of  \eqref{pin_vpgcf} is nonnegative, as shown in the proof of Theorem \ref{pinching}. By the maximum principle, we deduce that 
$$
\ln \frak d(t) \leq  \ln \frak d(0) - B \! \int_0^t  H_{min}(\tau) d\tau \to -\infty,
$$
where we have used Lemma \ref{0.19}. This allows us to conclude that $ \lim_{t \to \infty} \frak d(t) = 0,$
from which the proposition follows.
\edem

The above result guarantees that the limiting hypersurface, if it exists, has to be umbilical everywhere, and therefore is a round sphere. Accordingly, our remaining task consists of showing that a smooth limit actually exists, by proving that the principal curvatures cannot become zero in the limit.

\subsection{Interior H\"older estimate for the m$^{\text{\bf th}}$ mean curvature}

In order to prove our result, an essential step is the derivation of some kind of estimate on the curvature  (e.g. a Harnack inequality, or a H\"older estimate) which is uniform in time. Such an estimate allows us to say that, if the curvature is positive at a given point of our hypersurface, then it also satisfies a uniform lower bound in a whole neighborhood. However, there is a difficulty in deriving this type of inequalities, which has been pointed out in \cite{Sch2} and is related to the fact that the speed we are considering has a homogeneity degree larger than one in the curvatures; namely, we cannot ensure a priori that the evolution equations for the curvatures are uniformly parabolic. In fact, the operators $\Delta_c$ and $\Delta_{\dot \sigma}$ which appear in the equations, in contrast with the standard laplacian $\Delta$, become degenerate if the curvatures go to zero, and this is exactly the behavior we are not able to exclude at this stage.

Consequently, we will make use of the regularity theory for degenerate para\-bolic equations. Following
 \cite{Sch2}, we will prove a uniform $C^{\a}$-estimate for  the $m$-th mean curvature $H_m$ by means of Theorem \ref{tBF}, valid for equations of porous medium type. The procedure here is more complicate than in \cite{Sch2}; in particular, it is necessary to rewrite the evolution formula \eqref{Hm_pm} for  $H_m$ in a particular form (see the following lemma) which suits to the hypotheses of the regularity theorem.

\bl \lb{rew_ecK}
In a local coordinate system, the evolution equation for the $m^{\text{th}}$ mean curvature $H_m$ under \eqref{vpgcf} can be written as
\bec
\parcial{H_m}{t} = D_i\(\frac{\beta}{d} H_m^{\frac{1 - m}{m}} c^{ij} D_j H_m^d\) + \Gamma^j_{jl} c^{li} D_i H_m^\beta + (H_m^\beta - h) \tr_{c} (\a A), \lb{evK2}
\eec
where $d = \beta+\frac{m-1}{m}$ and $D_i$ denote derivatives  
with respect to the coordinates. 
\el

\bdem
First, the goal is to write the leading term of \eqref{Hm_pm} in a divergence form using the derivatives $D_i$. To do so, let us begin by writing the laplacian $\Delta_c$ in local coordinates: 
\begin{eqnarray}
\Delta_c & = & c^{ij} \left( D_i D_j  - \Gamma^l_{ij} D_l 
 \right) = D_i(c^{ij} D_j ) - (D_i c^{ij}) D_j - c^{ij} 
\Gamma^l_{ij} D_l \nn
\\& = & D_i(c^{ij} D_j) + \big(\Gamma^i_{il} c^{lj} + \Gamma^j_{il} c^{il} - \nabla_i c^{ij}\big) D_j  - c^{ij} 
\Gamma^l_{ij} D_l \nn
\\ & = & D_i(c^{ij} D_j) + \Gamma^i_{il} c^{lj} D_j, \lb{lap_cK}
\end{eqnarray}
where we have used part (e) in Lemma \ref{lpolS}.

On the other hand, if we set $d = \beta+\frac{m-1}{m}$, we have 
\bec  \lb{igu_der_pot}
D_j \sigma = \beta H_m^{\beta-1} D_j H_m = \frac{\beta}{d} H_m^{\frac{1 - m}{m}} D_j H_m^d.
\eec
Hence, by \eqref{lap_cK},
\begin{eqnarray*}
\Delta_c \, \sigma  = \frac{\beta}{d} D_i\(c^{ij} H^{\frac{1 - m}{m}} D_j H_m^d\)
+ \Gamma^i_{il} c^{lj} D_j \, \sigma.
\end{eqnarray*}
By substitution of this in \eqref{Hm_pm}, we reach formula \eqref{evK2} in the statement.
\edem

It is useful to estimate the tensor $c^{ij}$ introduced in  \eqref{def_c}. This is done conveniently in a coordinate system such that $g$ is the identity and $\alpha$ is diagonal at a given point. Then it is easy to see that
$$
c^{ij} = \frac{\partial H_m}{\partial k_i} \delta^{ij}.
$$
By definition of $H_m$, the derivative $\partial H_m/\partial k_i$ is a sum of products of $m -1$ principal curvatures. Hence we find, for any vector $\xi \in \re^n$,
$$
B'_{m,n} k_1^{m -1} |\xi|^2 \leq c^{ij} \xi_i \xi_j \leq B''_{m,n}  k_n^{m -1} |\xi|^2 
$$
for suitable constants $B'_{m,n},  B''_{m,n}$ depending only on $m,n$. It follows
$$
 B'_{m,n}(\eps\, H)^{m -1} |\xi|^2 \us{\eqref{i+ii}}{\leq} c^{ij} \xi_i \xi_j  \us{\eqref{kmenH}}{\leq} B''_{m,n}H^{m -1} |\xi|^2.
$$
For short, we express a double bound like the above one by writing $c^{ij} \approx H^{m-1} g^{ij}$. With this notation, we also obtain 
 \bec \lb{equi_key}
H_m^\frac{1-m}{m}c^{ij} \approx g^{ij},
\eec
because $H \approx H_m^{1/m}$ by part (d) of Lemma \ref{lpolS} and \eqref{HloF}. Another useful inequality related to $c$ is
\bec
\tr_c(\a A) = c^{ii} k_i^2 \leq B''_{m,n} H^{m - 1} \sum_i k_i^2 = B''_{m,n} H^{m-1} |A|^2 \leq B''_{m,n} C_1^{m + 1}, \lb{trc_aA}
\eec
which is true thanks to convexity and Proposition \ref{analogSch}.

We also need the following lemma.

\bl \lb{inte_kb} 
We can find a constant $\widehat C = \widehat C(n, m, \beta, M_0)$ satisfying
$$\int_{t_1}^{t_2} \int_M |\nabla H_m^d|^2 \, d\mu_t \, dt \leq \widehat C (1 + t_2 - t_1).$$
\el

\bdem
Let us begin by computing
\begin{align*}
\int_M |\nabla H_m^d|^2 \, d\mu_t  &
\ \ \usb{equi_key}{\approx} 
\int_M H_m^{\frac{m-1}{m}}   |\nabla H_m^d|_c^2 \, d\mu_t
\nn
 \\ & \ \! \us{\eqref{igu_der_pot}}{=}   \frac{d}{\beta}  \int_M \big<\nabla \sigma, \nabla H_m^d\big>_c  d\mu_t  \nn
= - \frac{d}{\beta}  \int_M H_m^d \, \Delta_c \sigma \, d\mu_t,
 \end{align*}
where the last equality follows using integration by parts and Lemma \ref{lpolS} (e). Next, we can use the evolution equation \eqref{Hm_pm} to deduce
\begin{align*}
\int_M |\nabla H_m^d|^2 \, d\mu_t  & \approx - \frac{d}{\beta(d +1)} \int_M \parcial{H_m^{d + 1}}{t}   \, d\mu_t + \frac{d}{\beta} \int_M (\sigma - h) H_m^d \, \tr_c (\a A)\, d\mu_t
\\ & \leq \mc C \deri{}{t} \int_M H_m^{d + 1} \, d\mu_t + \mc C'
\end{align*}
Notice that $\mc C'$ comes from \eqref{trc_aA}, \eqref{kmenH} and the bounds in Proposition \ref{analogSch}. Finally, recall that $0 < H_m \leq C$ and, by \eqref{db_inr}, $M_t$ is contained in a ball of radius $c_2$; these facts can be applied (after integrating the above inequality on $[t_1, t_2]$) to achieve the estimate in the statement.
\edem

\bp \lb{intHol_K}
For any point $(\bar x, \bar t) \in M \times (0, \infty)$, we can find a space-time neighborhood, say $\mc U \subset M \times (0, \infty)$, whose diameter does not depend on the point $(\bar x, \bar t)$ and such that 
$$ \big\|H_m\big\|_{C^\a(\mc U)} \leq \mc C,$$
for some constants $\mc C = \mc C(n, m, \beta, M_0) > 0$ and $\a \in (0, 1)$.
\ep

\bdem
We use the local parametrization of $M_t$ as the graph of a function $u:B_r \times J \to \re^{n+1}$ coming from Corollary \ref{graph_u}, where $J =\, \(\max\{\bar t -\tau, 0\}, \bar t + \tau\)$ for $\tau$ not depending on $\bar t$. It is not restrictive to assume that $\bar t>\tau$, since the $C^\alpha$ norm of $H_m$ on $M \times [0,\tau]$ is clearly finite by the compactness of $M$. 

If we consider $\tr_{c} (\a A)$ as a given function of $x \in B_r$,  \eqref{evK2} can be regarded as an equation  of the form \eqref{DiBF_ec} for  $v = H_{m}$, with 
\begin{align*}
 a^{ij} =\frac{\beta}{d} H_m^{\frac{1-m}{m}} c^{ij},  
  \quad \text{and } \quad f(x, t, v, D v) =  \Gamma^j_{jl} c^{li} D_i H_m^\beta + (H_m^\beta - h) \tr_c (\a A). 
\end{align*}

Next, notice that
 \bec 
 a(\xi, \xi) = \frac{\beta}{d} H_m^{\frac{1-m}{m}} c(\xi, \xi) \us{\eqref{equi_key}}{\approx} \frac{\beta}{d} |\xi|^2 \quad \text{ for any } \quad \xi \in \re^n, \lb{cota_a}
 \eec
and, setting $b_2: = 2 \, C_1^{m + 1} C^\beta$, a combination of \eqref{trc_aA} and Proposition \ref{analogSch} yields
 \begin{align}
 |f| &  \us{\eqref{igu_der_pot}}{\leq } \frac{\beta}{d}  \Big|
 \Gamma^j_{jl} c^{li} H_m^{\frac{1 -m}{m}} D_i H_m^d\Big| +  b_2 
  \us{\eqref{equi_key}}{\approx}   \frac{\beta}{d}|\Gamma_{jl}^j g^{li} D_i H_m^d| +  b_2  \nn 
 \ \leq  b_1 |D H_m^d| + b_2,
 \end{align}
where $b_1$ comes from  \eqref{invg_u} and the fact that $u$ is $C^2$-uniformly bounded (cf. Corollary \ref{graph_u}). Moreover, Lemma \ref{inte_kb} implies
$$\int \!\!\!\!\!\int_{B_r \times J} |D H_m^d|^2 \, d \mu_t \, dt \leq C(\tau).$$

Therefore, we are in position to apply Theorem \ref{tBF} with $r'=r/2$ and $\delta=\tau/2$ to deduce that $$\|H_m\|_{C^\a \left(B_{\frac r2}\times \left[ \, \bar t-  \frac{\tau}{2} \, , \,  \bar t + \frac{\tau}{2}\, \right] \right)} \leq \mc C$$
for suitable $\a \in \, (0, 1)$ and $0 < \mc C < \infty$ depending on $n, M_0, m, \beta$. \edem

\subsection{Exponential convergence to a round sphere}

\bt \lb{exp_kb}
Under the same hypotheses of Theorem \ref{lte_kb}, the $M_t$'s converge exponentially as $t \to \infty$  to a round sphere in the $C^\infty$-topology.
\et

\bdem
Let us take any sequence $\{\tau_j\} \subset [0, \infty)$ with $\tau_j \to \infty$. The uniform bounds on the curvatures imply that there exists a subsequence (again denoted by $\tau_j$) such that, up to translations,
$$X(\ccdot, \tau_j) \flecha X_\infty(\ccdot) \qquad \text{ in the } C^{1, \a}\text{-topology for any } \a < 1,$$
and $M_\infty:= X_\infty(M)$ is a convex $C^{1,1}$-hypersurface. In addition, by \eqref{db_inr}, at each time $\tau_j$ we can find a point $p_j \in M$ satisfying
$$H(p_j, \tau_j) \geq \frac{n}{c_2}.$$
Then \eqref{HloF} yields
\begin{equation}\label{estim0}
H_m(p_j, \tau_j) \geq \eps^m  H^m (p_j, \tau_j)  \geq   \(\frac{\eps \, n}{c_2}\)^m=: \gorro{C}(n, m, \beta, V) > 0
\end{equation}
for each fixed $j$. Proposition \ref{intHol_K} implies that $H_m$ cannot decrease too fast in the sense that  we can find a  $\delta > 0$ (independent of $(p_j, \tau_j)$) satisfying
\begin{equation} \label{estim}
H_m \big|_{B_\delta(p_j) \times [\tau_j -\delta,\tau_j+\delta]} \geq \frac{\gorro{C}}{2}.
\end{equation}

If $\delta$ is small enough, then $M_t \cap B_\delta(p_j)$ can be written as the graph of a function $u_j$ for any $t \in  [\tau_j -\delta,\tau_j+\delta]$ as in Corollary \ref{graph_u}. Now, we can repeat the arguments of the proof of Theorem \ref{teor_reg} on any neighborhood $B_\delta(p_j)  \times [\tau_j -\delta,\tau_j+\delta]$, but with upper bounds independent of time. In this way we obtain uniform $C^\infty$-estimates on the functions $u_j$ in suitable smaller neighborhoods, say of radius $\delta/2$. Therefore, we have that 
$$B_{\delta/2}(p_j) \cap M_{\tau_j} \flecha B_{\delta/2}(p_\infty) \cap M_\infty \quad \text{in } C^\infty,$$ 
where $X(p_j, \tau_j) \to p_\infty \in M_\infty$.

Recall that, by Proposition \ref{asym_pc}, the limit must be totally umbilic, and therefore is a portion of a sphere. By (\ref{estim0}), the sphere has $H_m$ curvature at least $\gorro C$. Then, in the neighbourhoods $B_{\delta/2}(p_j) \times [\tau_j -\delta/2,\tau_j+\delta/2]$ the $H_m$ curvature becomes arbitrarily close to a constant value not smaller than $\gorro C$. Using again the uniform H\"older continuity, we deduce that (\ref{estim}) holds, for $j$ large, in $B_{(3/2)\delta}(p_j)$ instead of $B_{\delta}(p_j)$. Thus we can extend the region where $M_\infty$ is known to be spherical. After a finite number of iterations, we deduce that $M_\infty$ is a sphere, whose radius is uniquely determined by the volume-preserving property.

Since the above argument can be applied to any sequence $\tau_j$, we conclude that the whole famlily $M_t$ converges to a sphere as $t \to \infty$, possibly up to a translation in space. This implies that $H_m$ tends to a positive constant as $t \to \infty$, and therefore it is bounded from below by some constant $\tilde \delta > 0$ for long times. Thus, we can repeat the argument in the proof of Proposition \ref{asym_pc} using
\begin{align*}
h \frac{q}{H}(|A|^2 - n H^2) &\geq B \, \tilde \delta \, (n^{-n} -q)\end{align*}
instead of \eqref{esti_th}. In this way we obtain that the rate of convergence of $q(\cdot,t)$ to $n^{-n}$ is exponential. With this property, we can argue exactly as in the proof of Theorem 3.5 in \cite{Sch2} to conclude that the second fundamental form of $M_t$ converges exponentially in $C^\infty$ to the one of a sphere. In particular, the $m$-th mean curvature satisfies an estimate of the form $|H_m^\beta(p,t)-h(t)| \leq {\cal C}e^{-{\cal C}'t}$ for all $p,t$. Then we have, for any $0<\tau_1 < \tau_2$, 
\begin{align*}
\max_M | X(x,\tau_2) - X(x,\tau_1)| & \le \max_M \int_{\tau_1}^{\tau_2} \left|\parcial{X}{t}(x,t)\right| dt  \\
&\leq \max_M \int_{\tau_1}^{\tau_2}  |h - H_m^{\beta}| \, dt \leq \frac{\cal C}{{\cal C}'} (e^{-{\cal C}'\tau_1}-
e^{-{\cal C}'\tau_2}), 
\end{align*}
which shows, by the Cauchy criterion, that $X(\cdot,t)$ has a limit as $t \to \infty$. Accordingly, the immersions converge to a round sphere with no need to add a translation. The $C^\infty$ convergence of the second fundamental form implies the $C^\infty$ convergence of the immersions and of the metric by standard arguments (see \cite{Sch2,An94}). \edem \bigskip

{\bf Acknowledgments:}  The first author was partially supported by the DGI (Spain) and FEDER Project MTM2007-65852, and by the net REAG MTM2008-01013-E. The second author was partially supported by the PRIN 2007WECYEA Project of the MIUR (Italy).

\bibliographystyle{amsplain}

\end{document}